\documentclass[a4paper]{article}
\usepackage{latexsym}
\usepackage{amssymb}
\usepackage{amsmath}
\newtheorem{conj}{Conjecture}[section]
\newtheorem{lemma}[conj]{Lemma}
\newtheorem{coro}[conj]{Corollary} 
\newtheorem{defi}[conj]{Definition}

\newtheorem{thm}{Theorem}

\newcommand{\qed}{\raisebox{-.8ex}{$\Box$}}
\newenvironment{bew}
{\noindent{\bf Proof.}}
{\hfill \qed}

\newcommand{\Aut}{{\rm Aut}}
\newcommand{\Out}{{\rm Out}}

\newcommand{\Syl}{{\rm Syl}}
\newcommand{\PGL}{{\rm PGL}}
\newcommand{\PSL}{{\rm PSL}}
\newcommand{\PSp}{{\rm PSp}}
\newcommand{\PGU}{{\rm PGU}}
\newcommand{\PSU}{{\rm PSU}}
\newcommand{\POmega}{{\rm P}\Omega}
\newcommand{\GL}{{\rm GL}}
\newcommand{\SL}{{\rm SL}}
\newcommand{\SU}{{\rm SU}}
\newcommand{\GU}{{\rm GU}}
\newcommand{\Sp}{{\rm Sp}}
\newcommand{\GF}{{\rm GF}}

\newcommand{\PGaL}{{\rm P} \Gamma {\rm L}}

\newcommand{\InnDiag} {~{\rm InnDiag}}
\newcommand{\Alt}{~{\rm Alt}}

\newcommand{\ZZ}{\mathbb{Z}}

\newcommand{\QR}{{\rm I}\kern-5.0pt {\rm Q} \kern2pt}

\newcommand{\refoverlineK}{\cite{Part1}(2.25)}
\newcommand{\refHKsuper}{\cite{Part1}(2.4)}

\newcommand{\refEnvelopeGroup}{\cite{Part1}~Thm~1~}
\newcommand{\refsolublegroups}{\cite{Part1}~Thm~2~}
\newcommand{\refHsuper} {\cite{Part1}(3.2)}
\newcommand{\refHisTGiffGisTG} {\cite{Part1}(3.6)}
\newcommand{\refOTTprime}{\cite{Part1}(3.7)}
\newcommand{\refTpowerorder}{\cite{Part1}(3.9)}
\newcommand{\refsubloops}{\cite{Part1}(3.12)}
\newcommand{\refnoHoverlineinvert}{\cite{Part1}(3.14)}

\newcommand{\refHeissPrime} {\cite{Part1}(3.17)}
\newcommand{\refZeroComponentCase}{\cite{Part1}(3.18)}
\newcommand{\refOupper}{\cite{Part1}(3.20)}
\newcommand{\refpassivecentralizingcomponents}{\cite{Part1}(3.21)}
\newcommand{\refTNloopEnvelope} {\cite{Part1}(3.24)}
\newcommand{\refNloopEmbedding}{\cite{Part1}(3.25)}
\newcommand{\refpcycles} {\cite{CommGraph}(4.1)}
\newcommand{\refEconnected} {\cite{CommGraph}~Thm~2~}
\newcommand{\refFconnected} {\cite{CommGraph}~Thm~2~}
\newcommand{\refsmallCC} {\cite{CommGraph}~Thm~4~}
\newcommand{\refalternatingcentralizers} {\cite{CommGraph}~(4.1)}

\title{On Bruck Loops of 2-power Exponent, II\thanks{This research is part of the project ``Transversals in Groups with an application to loops'' GZ: BA 2200/2-2 funded by
the DFG}}
\author{Alexander Stein}

\begin{document}

\maketitle

\begin{abstract}
As anounced in \cite{Part1}, we show that the non-passive finite simple groups are among
the $\PSL_2(q)$ with $q-1 \ge 4$ a 2-power.
\end{abstract}

\section{Introduction}
In part I, \cite{Part1},  the problem of {\em nice} loop folders to Bruck loops of 2-power exponent was reduced to a question
on simple groups. Specifically the problem arose, for which finite simple groups $S$ there is a BX2P-folder $(G,H,K)$
with $F^\ast(G/O_2(G)) \cong S$ and $G \ne O_2(G) H$.  We will give a partial answer here by using the classification of finite simple groups and its maximal subgroups, together with the results of \cite{CommGraph} and \cite{Part1}, see Theorem \ref{MainClassificationTheorem}:\\

Any such $S$ is among the groups $\PSL_2(q)$ for $q-1 \ge 4$ a 2-power.\\

It is known, that $\PSL_2(5)$ is such a group, but the question for other values of $q$ is open. The complexity of our proof is owed
to this situation: If for instance we knew that $\PSL_2(9)=\Alt_6$ is passive, our proof would shrink enormously. 
However we see no shortcuts if $\Alt_6$ were known to be non-passive.

In addition to the notation and definitions from \cite{Part1}, we use $\overline{G} :=G/O_2(G)$ and an 
overline denotes the image under the natural homomorphism from $G$ into $G/O_2(G)$.

The general setting of this paper is as follows:
$S$ is a finite simple group, $G$ a finite group with $F^\ast(G/O_2(G)) \cong S$. We look for $H \le G,K \subseteq G$, such that $(G,H,K)$ is a BX2P-folder. 
The overall strategy is to show, that $\overline{G}= \overline{H}$ by combining knowledge on loop folders with knowledge on the structure of
finite simple groups.          

\section{Some results on simple groups}
\begin{lemma}
\label{sporadiccyclicnormalizers}
Let $S$ be a sporadic simple group and $x\in S$, $o(x)=p$, $2<p$ prime. Then $|N_S(\langle x \rangle):C_S(x)|$ is not a 2-power, unless
$p$ is a Fermat prime. Thus in the non-Fermat case there exists an odd prime $s\mid p-1$ with  $s \mid |N_S(\langle x \rangle):C_S(x)|$. 
\end{lemma}

\begin{bew}
This lemma can be easily verified using the character tables in \cite{ATLAS}. 
The index $|N_S(\langle x \rangle):C_S(x)|$ determines the number $n_x$ of conjugacy classes of elements of order $p$ 
in $\langle x \rangle$. Recall, that $n_x= \frac{p-1}{|N_S(\langle x \rangle):C_S(x)|}$ and can be read off from 
the character tables, as the corresponding conjugacy classes have the same size and are powers of each other. )
\end{bew}

\subsection{Special centralizers}
We consider centralizers of elements of order 3 and 5 in the classical groups over $\GF(2)$. 

\begin{lemma}
\label{good35base}
Let $G \cong \SL_n(2)$, $\Sp_{n}(2)$, $\Omega^\pm_{n}(2)$ for $n \ge 2$,
$x \in G$, $o(x)=3$ or $5$,  
$V$ the natural $n$-dimensional $\GF(2)$ module for $G$.

Then $V = U_0 \oplus U_1 \oplus ... \oplus U_k$ with $U_0=C_V(x)$ and $U_i$ irreducible for $i>0$.
Moreover in the symplectic and orthogonal case, the direct summands can be choosen in such a way,
that $U_i \perp U_j $ for $i \ne j$ and the $U_i$ are nondegenerate. 
\end{lemma}

\begin{bew}
By coprime action the module splits into a direct sum of irreducibles. So suppose we have a nontrivial symplectic or quadratic form. 
We use induction on $\dim V$. By coprime action we have $[V,x] \perp C_V(x)=U_0$, 
so $C_V(x)=0$ by induction. 

Now $o(x)$ determines the minimal polynomial of $x$ uniquely:
It is $x^2+x+1$ for $o(x)=3$  and $x^4+x^3+x^2+x+1$ for $o(x)=5$ as irreducible $\GF(2)\langle x \rangle$-modules have dimension  2 resp. 4.

Let $U$ be some irreducible $x$-submodule of $V=[V,x]$. 
As $U$ is irreducible, either $U \cap U^\perp=0$ and $U$ is nondegenerate or
$U \cap U^\perp = U$, so $U$ is totally singular. 
Now $U^\perp/U$ is a $x$-module, but the extension splits over $U$, as $x$ acts semisimple.
So there exists an  $x$-invariant complement $W \le U^\perp$, which is nondegenerate as $U^\perp = U \perp W$.
By induction $W=0$, else we can produce the $U_i$ from proper subspaces $W$ and $W^\perp$. 

Therefore  $\dim V=4$ for $o(x)=3$ and $\dim V=8$ for $o(x)=5$. 
By inspection of the groups $\Omega^\pm_4(2)$, $\Sp_4(2)$, $\Omega^\pm_8(2)$ and $\Sp_8(2)$, see \cite{ATLAS},
these groups contain at most one class of fixed point free elements of order 3 resp.
5, except in case of $\Omega^+_8(2)$. In this case, there are three classes
of elements of order 5, which are transitively permuted by $\Out(G) \cong \Sigma_3$. 
In particular there is one class of elements of order 5, with $C_V(x) \ne 0$ and two fixed point free classes.
One of them comes from the $O_4^-(2) \perp O_4^-(2)$-decomposition, so for this element we have the above decomposition.
But the other class is an image under the graph automorphism of order 2, which preserves the module and quadratic form of $G$,
so we get a decomposition in this case too. 
Therefore there are irreducible and nondegenerate subspaces $U_1, U_2 \le V$ with $V = U_1 \perp U_2$.
\end{bew}

We get the following corollaries from basic representation theory:

\begin{coro}
\label{GL_n_2-35Centralizer}
Let $G \cong \GL_n(2)$ and $x,y \in G$ with $o(x)=3,o(y)=5$, $m=\dim [V,x], k=\dim [V,y]$.
Then $C_G(x) \cong \GL_{m/2}(4) \times \GL_{n-m}(2)$ and $C_G(y) \cong \GL_{k/4}(16) \times \GL_{n-k}(2)$. 
\end{coro}

\begin{coro}
\label{Sp_n_2-35Centralizer}
Let $G \cong \Sp_n(2)$ and $x,y \in G$ with $o(x)=3,o(y)=5$, $m=\dim [V,x], k=\dim [V,y]$.
Then $C_G(x) \cong \GU_{m/2}(2) \times \Sp_{n-m}(2)$ and $C_G(y) \cong \GU_{k/4}(4) \times \Sp_{n-k}(2)$. 
\end{coro}

In case of the orthogonal groups we formulate a weaker statement to avoid difficulties with automorphisms.

\begin{coro}
\label{OP_n_2-35Centralizer}
Let $G \cong \Omega^+_n(2)$ and $x,y \in G$ with $o(x)=3,o(y)=5$, $m=\dim [V,x], k=\dim [V,y]$.
Then $O^2(C_G(x)) \cong (\GU_{m/2}(2))' \times \Omega^{\varepsilon_1}_{n-m}(2)$ and 
$O^2(C_G(y)) \cong \GU_{k/4}(4) \times \Omega^{\varepsilon_2}_{n-k}(2)$
with $\varepsilon_1 = (-1)^{m/2}$ and $\varepsilon_2 = (-1)^{k/4}$.  
\end{coro}

\begin{coro}
\label{OM_n_2-35Centralizer}
Let $G \cong \Omega^-_n(2)$ and $x,y \in G$ with $o(x)=3,o(y)=5$, $m=\dim [V,x], k=\dim [V,y]$.
Then $O^2(C_G(x)) \cong (\GU_{m/2}(2))' \times \Omega^{\varepsilon_1}_{n-m}(2)$ and $O^2(C_G(y)) \cong \GU_{k/4}(4) \times \Omega^{\varepsilon_2}_{n-k}(2)$
with $\varepsilon_1 = (-1)^{1+m/2}$ and $\varepsilon_2 = (-1)^{1+k/4}$.  
\end{coro}

We now consider elements of order 3 and 5 in the groups $\PGU_n(2)$.  
Let $V$ be the natural $n$-dimensional $\GF(4)$-module of $\GU_n(2)$ and $\omega \in \GF(4)$ with $\omega^2+\omega+1=0$.

\begin{lemma}
\label{good_unitary_3_base}
Let $x \in \GU_n(2)$, such that $\overline{x} \in \PGU_n(2)$ has order 3. Then either 
\begin{itemize}
\item[(i)] $o(x)=3$ and  $x$ is diagonalizable. 
The eigenspaces to $1$, $\omega$ and $\omega^2$ are nondegenerate. 
\item[(ii)] $o(x)=9$, $x^3 \in Z(\GU_n(2))$ and $x$ has minimal polynomial $x^3-\omega$ or $x^3+\omega$.
There exist $x$-invariant subspaces $U_1,...,U_k$, $\dim U_i=3$ and $V = U_1 \perp U_2 ... \perp U_k$, $n=3k$ and the $U_i$ are nondegenerate.
\end{itemize}
\end{lemma}

\begin{bew}
Consider case (i) and let $u,v \in V$ eigenvectors to different eigenvalues $\lambda, \mu$.
Then $(u,v)=(u,v)^x = (u^x,v^x) = (\lambda u,\mu v) = \lambda \overline{\mu} (u,v)$. If $\lambda \ne \mu$ and
$\lambda, \mu \in \{ 1,\omega,\omega^2 \}$ this implies $(u,v)=0$, so the eigenspaces to different 
eigenvalues are orthogonal. As $x$ is diagonalizable, $V$ is the sum of these eigenspaces, so each
one is nondegenerate. 
Consider now case (ii). As $1 \ne x^3 \in Z(\GU_n(2))$, there are only the two choices $x^3= \omega Id$ or $x^3 = w^2 Id$
with $Id$ the identity matrix. Therefore the minimal polynomial is one of the two choices. As it is irreducible
we have $n$ a multiple of $3$, so $\SU_n(2)$ already contains $Z(\GU_n(2))$. Notice, that in this case
there are such elements, which come from the embedding $\GU_{n/3}(2^3).3 \le \GU_n(2)$, so there is a conjugacy class
of elements, which satisfies (ii). We show, that this subspace decomposition exists in general, by induction over $n$: 
Let $U$ be some irreducible $x$-submodule of $V$. Then either $U \cap U^\perp =0$ or $U \le U^\perp$, 
as $U$ is irreducible. If $U \cap U^\perp=0$, we proceed by induction on $U^\perp$. 
If $\dim U^\perp >3$, $U^ \perp$ has an $x$-invariant complement $W$ to $U$, as $x$ acts semisimple.
Then we proceed by induction on both $W$ and $W^\perp$, as $W \cap W^\perp=0$ and $V=W \perp W^\perp$. 
So $\dim V=6$. By \cite{ATLAS} there exists a unique conjugacy class of elements of order 9 in $\GU_6(2)$ 
with the property $x^3=\omega \in Z(\GU_6(2))$. It is class $3G$. By uniqueness it is the class, which allows
a decomposition into an orthogonal sum of two nondegenerate $x$-invariant subspaces, which comes from the embedding of $\GU_2(8) < \GU_6(2)$ , 
so the statement is proven.
\end{bew}

\begin{coro}
\label{unitary_3_centralizers}
Let $x \in G = \GU_n(2)$ with $x^3 \in Z(\GU_n(2))$. Then one of the following holds.
\begin{itemize}
\item[(i)] $o(x)=3$. Then $C_G(x) \cong \GU_{n_1}(2) \times \GU_{n_2}(2) \times \GU_{n_3}(2)$ with $n= n_1 + n_2 + n_3$. 
\item[(ii)] $o(x)=9$. Then $C_G(x) \cong \GU_{n/3}(8)$. 
\end{itemize}
\end{coro}

Notice that $\frac{t^5+1}{t+1}=t^4+t^3+t^2+t+1=(t^2+ \omega t +1)(t^2+\omega^2t+1)$ in $\GF(4)[t]$ with $\omega\in \GF(4)$, $\omega^2+\omega+1=0$.

\begin{lemma}
\label{good_unitary_5_base}
Let $x \in G=\GU_n(2)$ with $o(x)=5$. Then there exist $x$-invariant subspaces
$U, X_i,Y_i, i \in \{1..k\}$, $n=\dim U + 4 k$, with 
\begin{itemize}
\item  $U=C_V(x)$
\item  $\dim X_i= 2 =\dim Y_i$,
\item  $x$ is irreducible on $X_i$ and $Y_i$,
\item  $x$ has on $X_i$ the minimal polynomial $x^2+\omega  x +1$,
\item  $x$ has on $Y_i$ the minimal polynomial $x^2+\omega^2  x +1$, 
\item  $X_i \le X_i^ \perp$ and $Y_i \le Y_i^\perp$,
\item  $(X_i \oplus Y_i) \cap (X_i \oplus Y_i)^\perp=0$,
\item  $V = U \perp (X_1 \oplus Y_1) \perp ... \perp (X_k \oplus Y_k) $  
\end{itemize}
\end{lemma}

\begin{bew}
The proof proceeds by induction on $\dim V$. So we may assume, that $\dim U=\dim C_V(x)=0$. 
Let $X$ be some irreducible $x$-submodule, so $x$ is 2-dimensional. Then $X \le X^\perp$, 
as otherwise $X \cap X^\perp=0$, but $|\GU_2(2)|$ is not divisible by 5.
Let $W$ be an $x$-invariant complement to $X$ in $X^\perp$. Then $W$ is nondegenerate,
so if $W \ne 0$ the result holds by induction on both $W$ and $W^\perp$, so on $V$. 
If $W=0$, then $G=\GU_4(2)$ and $|\GU_4(2)|_5=5$. An easy calculation shows the statement in this case.
\end{bew}

\begin{coro}
\label{unitary_5_centralizers}
Let $x \in G = \GU_n(2)$, $o(x)=5$. Then $C_G(x) \cong \GU_k(2) \times \GU_{(n-k)/4}(4)$
with $k=\dim C_V(x)$. 
\end{coro}

\begin{bew}
This is a consequence of \ref{good_unitary_5_base}.
\end{bew}

\section{Passive simple groups: general arguments}
Here we give some general arguments, which we derive by combining knowledge about the finite simple groups with knowledge on BX2P-folders. 
These arguments are used in the next section to show, that almost all finite simple groups are passive. 

In this section $(G,H,K)$ is a BX2P-folder, $\overline{G}=G/O_2(G)$, $F^\ast(\overline{G}) \cong S$ with $S$ some finite simple nonabelian group, 
$S \le T \le \Aut(S)$ with $\overline{G} \cong T$ and $G_0$ the preimage of $F^\ast(\overline{G})$ in $G$. 

\subsection{An assumption and consequences}  

\begin{lemma}
\label{generated_by_involutions}
We may assume $G = \langle K \rangle$, so both $\overline{G} \cong T$ and $T/S$ are generated by involutions by \refoverlineK.
\end{lemma} 

\begin{bew}
This is \refHKsuper:
Let $g \in \langle K \rangle$. Then there exist $h \in H, k \in K$ with $g = h k$. As $k \in \langle K \rangle$,
$h \in \langle K \rangle \cap H$. Therefore $\langle K \rangle = K H_0$ with $H_0=\langle K \rangle \cap H$,
so $(\langle K \rangle,H \cap \langle K \rangle, K) $  is a subfolder of $(G,H,K)$.
As $|G:H|=|K|=|\langle K \rangle :H_0|$, this is a subfolder to the same loop. 
In particular $G \ne O_2(G) H$ iff $\langle K \rangle \ne O_2(\langle K \rangle)(H \cap \langle K \rangle)$. 
\end{bew}

This has consequences on the structure of $T$ and $T/S$: 

Recall that $|\Out(S)| \le 2$ for $S$ a sporadic or alternating group (different from $\Alt_6$.)
If $S$ is of Lie type, $\Out(S)$ can be more complicated, so we consider this case here.
By the famous Theorem of Steinberg on the structure of $\Aut(S)$, (Theorem 2.5.1 in \cite{GLS3}),
every automorphism of $S$ is a product of an inner, diagonal, field and graph automorphism. 
Moreover Theorem 2.5.12 in \cite{GLS3} gives a detailed description of $\Aut(S):$ 
$\Aut(S)$ is a semidirect product of a normal subgroup $\InnDiag(S) \le \Aut(S)$ with a subgroup $\Phi \Gamma$.
$\InnDiag(S)$ is the subgroup consisting of inner and diagonal automorphisms, while
$\Phi \Gamma$, is a product of a cyclic group $\Phi$ (inducing field automorphisms)
with a supplement $\Gamma$, such that $\Phi \Gamma /\Phi$ is a group of automorphisms
of the Dynkin diagram. 
By Theorem 2.5.12(e), if the group $S$ is untwisted and the Dynkin diagram contains only roots of one length,
then $\Phi \Gamma = \Phi \times \Gamma$ with $\Gamma$ the full automorphism group of the Dynkin diagram.
If the group is untwisted, but the Dynkin diagram contains roots of different length and a graph automorphism of order 2,
(so the group is $B_2(q), F_4(q)$ in characteristic 2 or $G_2(q)$ in characteristic 3) then 
$\Phi \Gamma$ is cyclic, with a generator in $\Gamma$, which squares to a Frobenius automorphism of $\GF(q)$
generating $\Phi$. 
If the group is twisted, then $\Gamma=1$. 
We will use definition 2.5.13 of \cite{GLS3} for the terms {\bf field}, {\bf graph-field} and {\bf graph} automorphism.
 
\begin{lemma}
\label{Out_restriction}
Let $S$ be a group of Lie type in characteristic $p$. 
If $T$ is generated by involutions, then $T/(T \cap \InnDiag(S))$ is isomorphic to $1, \ZZ_2, \ZZ_2 \times \ZZ_2$, $\Sigma_3$ or $\ZZ_2 \times \Sigma_3$.
By \ref{generated_by_involutions} we may assume this. In particular:
\begin{enumerate}
\item[(1)] $T$ does not contain field automorphisms of order bigger than two
\item[(2)] In case $S \cong B_2(q)$, $F_4(q)$, or $G_2(q)$, $|T:T \cap \InnDiag(S)|\le 2$.
\item[(3)] $|T:S|_2 \le 4$, if $S$ is a group of Lie type in characteristic $2$ or $|\InnDiag(S):S|$ is odd.
\end{enumerate}
\end{lemma} 

\begin{bew}
This is a consequence of Theorem of 2.5.12 of \cite{GLS3}.
\end{bew} 

We now establish some consequences in even characteristic:

\begin{lemma}
\label{Outer_O2}
Let $S$ be a group of Lie type in characteristic $2$ and $X \le T$ with $O_2(X) \cap S = 1$. 
Then $|O_2(X)| \le 4$,$O_2(X) \le Z(X)$ and $O_2(X) \cap O^2(X)=1=O_2(X) \cap O^2(T)$.  
\end{lemma}

\begin{bew}
By \ref{Out_restriction} the structure of $T/S$ is restricted. In particular we see, that $|O^2(T):S|$ is odd,
while $|O_2(X)S:S|=|O_2(X)|$. Therefore $|O^2(T) O_2(X)| = |O^2(T)| | O_2(X)|$,
so $O^2(T) \cap O_2(X)=1$. From \ref{Out_restriction}(3) we get $|O_2(X)| \le 4$.  
We can write $X = O^2(X) P$ for some $P \in \Syl_2(X)$. Then $|PO^2(T):O^2(T)| \le 4$, so $(PO^2(T))' \le O^2(T)$. 
As $X \le O^2(T)P$ we get $[X,O_2(X)] \le (PO^2(T))' \cap O_2(X) \le O^2(T) \cap O_2(X) = 1$. 
\end{bew}

\begin{coro}
\label{soluble_subloops_in_even_char}
Let $S$ be group of Lie type in characteristic $2$ and $U \le G$ a subgroup to a soluble subloop. 
If $\overline{U}$ is reductive, so $O_2(\overline{U}) \cap \overline{G_0}=1$, then $O^2(\overline{U}) \le \overline{H}$. 
\end{coro}

\begin{bew}
By \ref{Outer_O2} we get $O_2(U) \cap O^2(U) \le O_2(G)$. Now \refOupper gives the statement.
\end{bew}

Standard examples, where \ref{soluble_subloops_in_even_char} can be applied, are centralizers of elements of odd
order in $\overline{H} \cap \overline{G_0}$, as centralizers of semisimple elements in $S$ are reductive. 

\subsection{Centralizers of odd order elements}
The main connection between the local structure of loops and local subgroups of almost simple groups is \refsubloops(i).
Given a subgroup $1 \ne L \le H$ of odd order, we have $C_G(L)$ covered by $C_{\overline{G}}(\overline{L})$ due to coprime action. 
As $C_G(L)$ does not cover $\overline{G}$, we can apply \refEnvelopeGroup on $\langle C_G(L) \cap K \rangle$. On the other hand
we know from the local structure of simple groups, how $C_{\overline{G}}(\overline{L})$ looks like. Putting things together,
we often can identify $C_{\overline{H}}(\overline{L})$ within $C_{\overline{G}}(\overline{L})$, without even knowing $\overline{H}$ completely.
We give here some lemmata, based on this idea. 

\begin{lemma}
\label{alternating_centralizers}
Let $S\cong \Alt_n$ for $n \ge 7$ , $x \in H$ be of odd prime order $p$ and $k$ be the number of fixed points of 
$\overline{x}$ in the natural action of $\overline{G}$ on $n$ points. 
\begin{itemize}
\item If $k \ne 5$, then $|C_G(x):C_H(x)|$ is a 2-power. 
\item If $k \not\in \{ 4,5 \}$, then $O^2(C_{\overline{G}}(\overline{x})) \le \overline{H}$. 
Remember that $O_p(C_{\overline{G}}(\overline{x})) \le \overline{H}$ by \refOTTprime.
\item If $k=5$, then $|C_G(x)O_2(G):C_H(x)O_2(G)| \in \{1,6 \}$. 
\item In any case, $\overline{H}$ contains $p$-cycles from $O_p(C_{\overline{G}}(\overline{x}))$. 
\end{itemize}
\end{lemma}
 
\begin{bew}
The structure of $C_{\overline{G}}(\overline{x})$ is well known, as elements commuting with $\overline{x}$
permute the cycles of $\overline{x}$ and act on the $k$ fixed points, see also \refalternatingcentralizers. 

So we apply the structure description of \refEnvelopeGroup on $C_G(x)$. Now $C_G(x)$ is an extension of $O_2(G) \cap C_G(x)$
by $C_{\overline{G}}(\overline{x})$. In particular $C_G(x)/O_2(C_G(x))$ has no subnormal $\PGL_2(q)$ for $q>5$. There may or may not be a
subnormal $\PSL_2(9) \cong \Alt_6$, but the outer involution is missing in $C_{\overline{G}}(\overline{x})$, seen as a subgroup of $\Sigma_n$.
Therefore the subloop to $C_G(x)$ is soluble, if $k\ne 5$, as there is no subnormal $\Sigma_5$ or $\PGL_2(9)$ in this case.

If $k =5$, there may be a subnormal $\Sigma_5$, acting on the 5 fixed points of $\overline{x}$. In this case
\refEnvelopeGroup describes the structure of $\langle \overline{K} \cap C_{\overline{G}}(\overline{x})$ and
its intersection with $\overline{H}$. 

If $k \ne 4,5$, then $O_2(C_{\overline{G}}(\overline{x}))=1$ and the subloop to $C_G(x)$ is soluble, so by \refOupper
we have $O^2(C_{\overline{G}}(\overline{x})) \le \overline{H}$.

We could even determine, which elements end up in $\overline{H}$ in cases $k=4$ and $k=5$, but have no use for it. 

By \refOTTprime, $\overline{H}$ contains $O_p(C_{\overline{G}}(\overline{x}))$, so in particular $p$-cycles.       
\end{bew}

\begin{lemma}
\label{sporadic_centralizers}
Let $S$ be sporadic and $x \in H$ be of odd prime order $p$. Then $|C_G(x):C_H(x)|$ is not divisible by $p$,
unless maybe $(p,S)$ is one of $(3,M_{23})$, $(3,HS)$ or $(5,Suz)$.  

\end{lemma}

\begin{bew}
This is a consequence of the list of centralizers of elements of prime orders
in sporadic groups in \cite{GLS3} and \refEnvelopeGroup: 
In case $p=3$  we have to check, which centralizers of elements of order 3 contain components of type $\PSL_2(r)$ 
for $r$ some Fermat prime, $r \ge 5$. 
In case $p=5$ we have to check, that centralizers of elements of order 5 do not contain components of type $\Alt_6$. 
In the cases listed above, there are elements of order 3 resp. 5, such that the corresponding centralizers
may even contain subnormal subgroups isomorphic to $\PGL_2(5)$ resp. $\PGL_2(9)$, depending on the presence
of outer automorphisms of $S$ in $T$. 
\end{bew}

\begin{lemma}
\label{odd_char_centralizers}
Let $S$ be a group of Lie type in odd characteristic $p$ and $x \in H$ be of odd prime order $r$.
If $r=p$, then $|C_G(x):C_H(x)|$ is not divisible by $p$. 
If $r \ne p$, then $\overline{H} \cap \overline{G_0}$ contains elements of order $p$ or 
$C_{\overline{G}}(\overline{x})$ is soluble or both.    
\end{lemma}

\begin{bew}
By \ref{Out_restriction}, $\overline{x}$ is either innerdiagonal or $S \cong D_4(q)$ and $\overline{x}$ is a graph or graph-field
automorphism of order 3. 
In case $\overline{x}$ innerdiagonal we have the cases $\overline{x}$ unipotent ($r=p$) or semisimple ($r \ne p$). 

In the semisimple case we use Theorem 4.2.2 of \cite{GLS3} for the description of $C_{\overline{G}}(\overline{x})$.
In particular nonsoluble composition factors of $C_{\overline{G}}(\overline{x})$ are groups of Lie type in characteristic $p$.
Therefore, if $\langle C_G(x) \cap K \rangle \not\le O_2(C_G(x))$, then by \refEnvelopeGroup and induction,
$\overline{G_0} \cap \overline{H}$ contains elements of order $p$. Recall, that this happens only, if
$C_G(x)$ is nonsoluble by \refsolublegroups. 
  
If $\overline{x}$ is unipotent, so $r=p$ and $x \in G_0$, we apply the Borel-Tits Theorem on $C_{\overline{G}}(\overline{x})$,
(Theorem 3.1.3 and Corollary 3.1.4 of \cite{GLS3}) and \refZeroComponentCase.

In the exceptional case of $D_4(q)$, the centralizer of a graph or graph-field automorphism of order 3 
is described by Proposition 4.9.1 and 4.9.2 of \cite{GLS3} and listed in 4.7.3 in \cite{GLS3}. 
$C_{\overline{G}}(\overline{x})$ contains a subnormal ${}^3D_4(q^{1/3})$ or a $G_2(q)$, 
which by induction is passive, so the subloop to $C_G(x)$ is soluble, so $\overline{G_0} \cap \overline{H}$ contains
elements of order $p$.   
\end{bew}
 
\begin{lemma}
\label{even_char_centralizers}
Let $S$ be a group of Lie type in even characteristic with $q$ the field parameter of $S$. (If $S$ is defined relative to a field extension,
$q$ is the size of the smaller field.)
Let $x \in H$ be of odd prime order $r$. 
\begin{itemize}
\item If $q\ge 8$, then $C_G(x)$ gives a soluble subloop. 
\item If $q=4$, then $\pi(|C_G(x):C_H(x)|) \subseteq \{2,3 \}$, so the subloop to $C_G(x)$ may not be soluble.
\item If $q=2$, then $\pi(|C_G(x):C_H(x)|) \subseteq \{2,3,5\}$, so again the subloop to $C_G(x)$ may not be soluble.
\item If the subloop to $C_G(x)$ is soluble, then $O^2(C_{\overline{G}}(\overline{x})) \le \overline{H}$.  
\end{itemize}
\end{lemma}

\begin{bew}
By \ref{Out_restriction}, $T$ does not contain field automorphisms of odd order.
Therefore automorphisms of $S$ of odd order are either innerdiagonal, so semisimple or
are graph or graph-field automorphisms of order 3 in case of $S\cong D_4(q)$. 

In case  $S \cong D_4(q)$ and $\overline{x}$ induces a graph or graph field automorphism of order 3,
we refer to Propositions 4.9.1 and 4.9.2 as well as 4.7.3 of \cite{GLS3} for the centralizer of these automorphisms.
In particular $C_{\overline{G}}(\overline{x})$ is reductive and contains a unique component isomorphic to $G_2(q)$ or ${}^3D_4(q^{1/3})$. 
By induction this component is passive, so the subloop to $C_G(x)$ is soluble.  

If $\overline{x}$ induces a semisimple, innderdiagonal automorphism, its centralizer is reductive too. Moreover the structure of the centralizer is described by Theorem 4.2.2 of \cite{GLS3}.
In particular the nonsoluble composition factors come from components, which are groups of Lie type in characteristic 2 and
defined over field extensions of $\GF(q)$. The only nonpassive components, which may arise, are $\Sp_4(2)' \cong \Alt_6 \cong \PSL_2(9)$ and
$\PSL_2(4) \cong \Alt_5 \cong \PSL_2(5)$. As $\Sp_4(2)'$ is a group defined over $\GF(2)$, it does not arise if $q>2$. As $\PSL_2(4)$
is a group defined over $\GF(4)$, it does not arise as a component, if $q>4$. This is the reason for the case division $q>4$, $q=4$ and $q=2$. 
Now \refEnvelopeGroup gives solubility of the subloop for $q>4$, as no such component occurs.

In case $q=4$ it gives, that $|C_G(x):C_H(x)|$ is a 2-power times a 3-power,
as only $\PSL_2(4)$-components may be not passive.

Finally in case $q=2$ it gives, that $\pi(|C_G(x):C_H(x)|) \subseteq \{ 2,3,5 \}$,
as there may occure $\PSL_2(4)$- or $\Sp_4(2)'$-components, but no other non-passive components.

Notice, that in any case $C_{\overline{G}}(\overline{x})$ is reductive, so we can use \ref{soluble_subloops_in_even_char}, if the subloop is soluble.    
\end{bew}

\subsection{The property $FS_p$}
Recall the class ${\cal L}_S$ from {\cite{Part1} (4.1). We have to generalize
this concept slightly to our group $T$ to avoid difficulties: 
\begin{defi}
We denote with ${\ell}_T$ the class of Bruck loops of 2-power exponent,
for which a BX2P-folder $(G_X,H_X,K_X)$ exists with $G_X/O_2(G_X) \cong T$.   
Let $p \in \pi(T), p>2$. The class ${\ell}_T$ has the {\bf property} $FS_p$, iff for all $X \in {\ell}_T$:
Either $|X|_p = |T|_p$ or $|X|_p=1$. 
\end{defi}

\begin{lemma}
The class ${\ell}_T$ has property $FS_p$, iff for every BX2P-folder $(G,H,K)$ 
with $G/O_2(G) \cong T$: $p \nmid (|H|,|G:H|)$, so either $\Syl_p(H) \subseteq \Syl_p(G)$ or $p \nmid |H|$. 
(And $FS_p$ stands for 'full Sylow-$p$'.)  
\end{lemma}

\begin{bew}
Suppose $(G,H,K)$ is a BX2P-folder with the property: 
If $p \in \pi(H)$, then $p \nmid |G:H|$. 
Then either $p \in \pi(H)$ and $p \nmid |X|=|G:H|$ or $p \not\in \pi(H)$, so $|X|_p=|G:H|_p=|G|_p=|T|_p$.
If every BX2P-folder $(G,H,K)$ with $G/O_2(G) \cong T$ has the above property, 
then the class ${\ell}_T$ has property $FS_p$. 
The converse statement is immediate from the definition. Notice, that 
the property $(|H_X|,|G_X:H_X|)_p=1$ depends only on the isomorphism type of $X$, not on the 
particular BX2P-folder $(G_X,H_X,K_X)$ to $X$.
\end{bew}

The reason for defining this property $FS_p$ is, that it can be established
from the $p$-local structure of $T$ in many cases, and has powerful applications. 

\begin{lemma}
\label{FS_p_criterion}
The class ${\ell}_T$ has property $FS_p$, if for any BX2P-folder $(G,H,K)$ with $\overline{G}=G/O_2(G) \cong T$ and any $x \in H$ with $o(\overline{x})=p$ one of the following conditions is satisfied: 
\begin{enumerate}
\item[(0)] $p \nmid |C_{G}(x):C_{H}(x)|$. 
\item[(1)] $C_{\overline{G}}(\overline{x})$ is soluble. 
\item[(2)] $F^\ast(C_{\overline{G}}(\overline{x})) = O_p(C_{\overline{G}}(\overline{x}))$ for $p>2$. 
\item[(3)] $C_{\overline{G}}(\overline{x})/O_2(C_{\overline{G}}(\overline{x}))$ has only passive components.
\item[(4)] $C_{\overline{G}}(\overline{x})/O_2(C_{\overline{G}}(\overline{x}))$ has no subnormal $\PGL_2(q)$ for $q=9$ or a Fermat prime $q\ge 5$
with $p | q+1$.
\end{enumerate}
Notice, that conditions (1)-(4) depend only on the structure of $T$, not on a particular BX2P-folder.
\end{lemma}
\begin{bew}
The general argument in all cases is the same:
Suppose $p \in \pi(H)$. We will show, that each of (1)-(4) implies (0).
Once this is established, $H$ contains a Sylow-$p$-subgroup of $G$ for the following reason: 
Every element $x\in G$ of order $p$ is centralized by some element $y$ of order $p$ with 
$y \in \Omega_1(Z(Y))$ for some $Y \in \Syl_p(G)$ with $x \in Y$.
If $x \in H$, then by (0), some $G$-conjugate $z$ of $y$ is in $H$. 
Using (0) on $C_G(z)$ we get a Sylow-$p$-subgroup of $G$ into $H$.
Notice, that $C_T(\overline{x})$ is covered by $C_G(x)$ for $x\in H$ some preimage of $\overline{x}$ 
of order $p$, due to coprime action. This enables to establish (0) from
information of $\overline{G}$ only: 

By \refsubloops, $C_G(x)$ gives a subfolder, so we can use inductive arguments on $C_G(x)$. 
In case (1), if $C_{\overline{G}}(x)$ is soluble, then $C_G(x)$ is soluble, so by \refsolublegroups $|C_G(x):C_H(x)|$
is a 2-power and we have (0).
In case (2) we have $|C_G(x):C_H(x)|_{2'}=1$ by \refZeroComponentCase.
In case (3) we use \refEnvelopeGroup to establish,
that $|C_G(x):C_H(x)|$ is a 2-power, as only passive components show up, so 
$\langle K \cap C_G(x) \rangle \trianglelefteq C_G(x)$ has to be a 2-group by \refTpowerorder.
Case (4) gives the most powerful criterion: The condition on $C_{\overline{G}}(\overline{x})$
gives a condition on the structure of $C_G(x)$. Using the factorization $C_G(x) = C_H(x) \langle C_G(x) \cap K \rangle$ and 
\refEnvelopeGroup for the structure description of $\langle C_G(x) \cap K \rangle$, we see that 
$|O_2(G) C_G(x):O_2(G) C_H(x)|$ is a product of integers $q_i+1$ for $q_i-1\ge 4$ a 2-power.
But the condition (4) on $C_{\overline{G}}(\overline{x})$ ensures, that none of $q_i+1$ is divisible
by $p$, so $|C_G(x):C_H(x)|$ is indeed not divisible by $p$.
\end{bew}

We now establish the $FS_p$-property in certain cases of the classification of finite simple groups. 
Notice, that a critical point may arise from the existence of outer automorphisms of $S$
of odd order, as we may have to establish condition (0) for such elements too. This is one reason for the assumption,
which led to \ref{Out_restriction}. 

\begin{lemma}
\label{FS_p_alternating}
Let $S$ be an alternating group and $S \le T \le \Aut(S)$ and $p>3, p \in \pi(T)$.  
Then ${\ell}_T$ has property $FS_p$. 
\end{lemma}

\begin{bew}
By \ref{alternating_centralizers} we have condition (4) of \ref{FS_p_criterion} for $p>3$.
\end{bew}  

\begin{lemma}
\label{FS_p_sporadic}
If $S$ is sporadic, $S \le T \le \Aut(S)$ and $p>2$, then ${\ell}_T$ has property $FS_p$,
unless $(p,S)$ is one of $(3,M_{23})$, $(3,HS)$  or $(5,Suz)$. 
\end{lemma}

\begin{bew}
We can use condition (4) of \ref{FS_p_criterion} by \ref{sporadic_centralizers}.
\end{bew}

\begin{lemma}
\label{FS_characteristic_p}
If $S$ is a group of Lie type in characteristic $p$, $p>2$ and $S \le T \le \Aut(S)$,
then ${\ell}_T$ has property $FS_p$. 
\end{lemma}

\begin{bew}
By \ref{odd_char_centralizers} we have either condition (2) or condition (3) of \ref{FS_p_criterion}.
\end{bew} 

This fact implies later, that in odd characteristic $p$ and Lie rank at most two,
$\overline{G}=\overline{H}$, if $p \in \pi(H)$, see \ref{OddChar}. This implies then the $FS_r$-property for primes $r \ne p$.
But before this we continue with groups of Lie type in characteristic 2.

\begin{lemma}
\label{FS_p_even_char}
Let $S$ be a group of Lie type in characteristic $2$, 
defined over the field with $q$ elements. (In case the group is defined relative to a field extension,
$q$ refers to the smaller field.) Let $S \le T \le \Aut(S)$. 
If $q>4$, then ${\ell}_T$ has property $FS_p$ for every prime $p>2$. 
If $q=4$, then ${\ell}_T$ has property $FS_p$ for every prime $p>3$. 
If $q=2$, then ${\ell}_T$ has property $FS_p$ for every prime $p>5$. 
\end{lemma}

\begin{bew}
This is a consequence of \ref{even_char_centralizers}, which enables condition (4) of \ref{FS_p_criterion}
under the given restrictions.
\end{bew}

\begin{coro}
\label{commuting_in_char_2} 
Let $S$ be a group of Lie type in characteristic $2$, as in \ref{even_char_centralizers}.
For $q>4$ we have $O^2(C_{\overline{G}}(\overline{x})) \le \overline{H}$, so if $\overline{x},\overline{y} \in \overline{G}$
are elements of odd order with $[x,y]=1$ and $\overline{x} \in \overline{H}$, then $\overline{y} \in \overline{H}$.
For $q=4$ we have either $FS_3$-property or $\overline{H}$ contains elements of order 15.  
Furthermore either $5 \in \pi(H)$ or with $\overline{x} \in \overline{H}$ of odd prime order the full connected
component ${\cal C}_{\overline{x}}$ of $\overline{x}$ in $\Gamma_{\cal O}$ to $S$ is contained in $\overline{H}$. 
For $q=2$ we have either both $FS_3$ and $FS_5$-property or $\overline{H}$ contains elements of order 15.
\end{coro}

\begin{bew}
For $q>4$, subloops to centralizers of elements of odd prime order are soluble. 
Then the statement is \ref{soluble_subloops_in_even_char} together with the fact,
that centralizers of semisimple elements (or outer automorphisms of order 3 in case of $D_4(q)$)
are reductive, by Theorems 4.2.2, 4.9.1, 4.9.2  and 4.7.3 of \cite{GLS3}.

For $q=4$, how can $FS_3$-property fail? Only, if there is some element $x \in H$, $o(x)=3$,
such that $C_G(x)$ gives a nonsoluble subloop in $G$, so $C_{\overline{G}}(\overline{x})$
contains a subnormal $\PSL_2(4)$. In that case the size of the subloop is a 2-power times a 3-power,
so $C_H(x)$ contains elements of order 5, so $\overline{H}$ contains elements of order 15.
If $\overline{H}$ contains no elements of order 5, the subloops to $C_G(x)$ for $x \in H$ of odd order are 
soluble, so by \ref{even_char_centralizers} ${\cal C}_{\overline{x}}\subseteq \overline{H}$. 
For $q=2$ there is also the possibility for the $FS_5$-property to fail: There may exist
some element $x \in H$, $o(x)=5$, such that the subloop to $C_G(x)$ is nonsoluble, 
but the size of the loop is divisible by 5. So some elements of order 5 are commutators of 
elements of $\overline{K}$, and $\overline{H}$ cannot contain a Sylow-5-subgroup of $\overline{G}$.
In that case $C_{\overline{G}}(\overline{x})$ contains a subnormal $\Sp_4(2)'$, so $C_H(x)$ contains
elements of order 15.  
\end{bew}

\subsection{Terminal elements}
One strategy in the generic case (where simple groups are big enough), is the identification of terminal elements. 
Plainly, an element $\overline{x} \in \overline{G}$ is terminal, if $\overline{x} \in \overline{H}$ implies $\overline{G}=\overline{H}$. 
This property can sometimes established from the structure of $C_T(\overline{x})$ together with the structure of $T$. 
We will give here some examples, which we will use in the next section.
We need a little lemma:

\begin{lemma}
\label{simple_group_is_enough}
Assume $\overline{G_0} \le \overline{H}$. Then $\overline{G}=\overline{H}$.
\end{lemma}

\begin{bew}
Assume otherwise, so the loop to $G$ is nonsoluble. We get a contradiction from \refNloopEmbedding:
If the loop to $G$ is nonsoluble, then $\overline{G_0}$ contains elements of odd order, which are not in $\overline{H}$, 
as they are commutators of elements of $\overline{K}= \overline{\Lambda}$.   
\end{bew}

\begin{lemma}
\label{terminal_three_cycles}
Let $S \cong \Alt_n$ for $n \ge 9$. If $\overline{H}$ contains a 3-cycle $\overline{x}$, then $\overline{H}=\overline{G}$.  
\end{lemma}

\begin{bew}
Let $x \in H$ be of order 3, a preimage of $\overline{x}$. By \ref{alternating_centralizers} we have, that $C_G(x)$ gives
a soluble subloop. Moreover $O_2(C_{\overline{G}}(\overline{x}))=1$, so by \refOupper we have $O^2(C_{\overline{G}}(\overline{x})) \le \overline{H}$. 
In particular $\overline{H}$ contains with $\overline{x}$ all 3-cycles, which commute with $\overline{x}$. As the commuting graph of 3-cycles is connected (see \refpcycles), we have $\overline{G_0} \le \overline{H}$, which implies $\overline{G}=\overline{H}$ by \ref{simple_group_is_enough}.
\end{bew}

\begin{lemma}
\label{OddChar_I}
Let $S$ be a group of Lie type in odd characteristic $p$. Assume the (twisted) Lie rank is not 1
(so $S$ is not of type $A_1$, ${}^2A_2$ or ${}^2G_2$) and $\overline{H}$ contains elements of order $p$. 
Then $\overline{G}=\overline{H}$. 
\end{lemma}

\begin{bew}
By \ref{FS_characteristic_p}, $\ell_T$ has property $FS_p$. Since $\overline{H}$ contains elements of order $p$, 
we can choose $P \in \Syl_p(\overline{H}) \subseteq \Syl_p(\overline{G})$. As the (twisted) Lie rank of $S$ is not 1, 
we find subgroups $V_1,V_2 \le P$ with $P_i:=N_{\overline{G}}(V_i)$ and $\overline{G_0} = \langle P_1,P_2 \rangle$. 
By \refsubloops(1) and \refZeroComponentCase, $P_i \le \overline{H}$, so $\overline{G_0} \le \overline{H}$. 
By \ref{simple_group_is_enough} we have $\overline{G}=\overline{H}$.  
\end{bew}

\begin{coro}
\label{OddChar}
Let $S$ be a group of Lie type in odd characteristic $p$. 
Assume, the (twisted) Lie rank is not 1 and let $r \in \pi(T), r>2$. 
\begin{itemize}
\item[(i)]  $\ell_T$ has property $FS_r$.
\item[(ii)] Suppose $\overline{H} \ne \overline{G}$. Let $x \in H \cap G_0$ of odd prime order. 
Then $\overline{x}$ is not in the big connected component of $\Gamma_{\cal O}$ (provided, that component exists). 
\end{itemize}
\end{coro}

\begin{bew}
(i): Let $\overline{x} \in \overline{H}$ be of order $r$, $p \ne r \ne 2$. By \ref{odd_char_centralizers}
either $\overline{H} \cap \overline{G_0}$ contains elements of order $p$ or $C_{\overline{G}}(\overline{x})$
is soluble (or both). In the first case $\overline{H}=\overline{G}$ by \ref{OddChar_I},
so $\overline{H}$ contains a Sylow-$r$-subgroup and we have (0) of criterion \ref{FS_p_criterion}, 
while in the second case we use (1) of \ref{FS_p_criterion} to get property $FS_r$. 
The first case depends on the fact, that the groups $\PSL_2(q)$ for $q\ge 5$ a Fermat prime or $q=9$ cannot be seen as groups of Lie type in another odd characteristic.

(ii): Notice, that $H \cap G_0$ contains elements of odd order by \refNloopEmbedding. 
Unfortunately we cannot use \refOupper or \refpassivecentralizingcomponents, 
as we have no control on $O_2(C_{\overline{G}}(\overline{x}))$. 
Suppose $\overline{x}$ is in the big connected component of $\Gamma_{\cal O}$.
Let $\pi=(\overline{x_i})$, $i \in \{1,..,k\}$ be a path of shortest length in $\Gamma_{\cal O}$ 
from some element $\overline{x_1} \in \overline{H}$ of odd prime order to some element $\overline{x_k}$ of order $p$. 
Suppose $s=o(\overline{x_1})=o(\overline{x_2})$. 
By $FS_s$-property, $H$ contains a Sylow-$s$-subgroup of $G$, so we find a $\overline{g} \in \overline{G}$, such that
$\langle \overline{x_1},\overline{x_2} \rangle ^{\overline{g}}  \le \overline{H}$. As $\overline{x_2}^{\overline{g}} \in \overline{H}$,
we get a shorter path by dropping $\overline{x_1}^{\overline{g}}$ from $\pi^{\overline{g}}$. 
Suppose $s=o(\overline{x_1}) \ne o(\overline{x_2})=t$. Choose $x_1 \in H$ in the preimage of $\overline{x_1}$.
Recall, that the subloop to $C_G(x_1)$ is soluble.
Furthermore $C_G(x_1)$ covers $C_{\overline{G}}(\overline{x_1})$ by coprime action, so $C_G(x_1)$ contains elements of order $t$.
As $|C_G(x_1):C_H(x_1)|$ is a 2-power, $t \in \pi(H)$ and by $FS_t$-property, $H$ contains a Sylow-$t$-subgroup of $G$.
Therefore some $\overline{g} \in \overline{G}$ exists with $\overline{x_2}^{\overline{g}} \in \overline{H}$. 
We get again a shorter path from $\pi^{\overline{g}}$ by dropping $\overline{x_1}^{\overline{g}}$.
Consequently the path consists of $\overline{x_1}$ only, so $\overline{H}$ contains elements of order $p$
and $\overline{H}=\overline{G}$. 
\end{bew}

\begin{lemma}
\label{EvenCharTerminals}
Let $S$ be a group of Lie type in characteristic $2$ and $x \in H \cap G_0$ of odd order $r>1$. 
Assume, that the commuting graph of $\overline{x}^{\overline{G_0}}$ in $\overline{G_0}$ is connected and the subloop to $C_G(x)$ is soluble. 
Then $\overline{H}=\overline{G}$. 
\end{lemma}

\begin{bew}
By \ref{soluble_subloops_in_even_char} we have $O^2(C_{\overline{G}}(\overline{y})) \le \overline{H}$ for $y=x$ and conjugates of $x$, which are contained in $\overline{H}$. 
Therefore with $\overline{x} \in \overline{H}$ all $\overline{G_0}$-conjugates of $\overline{x}$, which commute with $\overline{x}$, are in $\overline{H}$ too.
Then $\overline{H}$ contains $\overline{G_0}$, so by \ref{simple_group_is_enough}, $\overline{H}=\overline{G}$. 
\end{bew}

Once certain elements are established as being terminal, we can classify 'isolated elements'. An element $x \in \overline{G}$ is called inductive,
if $\overline{x} \in \overline{H}$ implies, that $\overline{y} \in \overline{H}$ for $\overline{y}$ either a terminal element or an element already identified as being inductive. 
Elements, which are neither terminal nor inductive are called isolated.  In the odd characteristic case and the characteristic 2-case with $q>4$,
inductive elements are simply elements from the same connected component in $\Gamma_{\cal O}$, while isolated elements
come from small connected components. 
 
\subsection{Other recurring arguments}

The following lemma is often used in case of cyclic groups $\overline{L} \le \overline{H}$ to get additional primes into $\overline{H}$.
Therefore most small connected components of $\Gamma_{\cal O}$ give inductive elements. 

\begin{lemma}
\label{H_covers_automizer}
Given $1 \ne \overline{L} \le \overline{H}$ with $|\overline{L}|$ odd, then $|N_{\overline{G}}(\overline{L}):C_{\overline{G}}(\overline{L})|_{2'}$ divides $|\overline{H}|$. 
\end{lemma}

\begin{bew}
Let $L \le H$ be a preimage of $\overline{L}$ with $|L|= |\overline{L}|$. 
By \refsubloops(i), both $(N_G(L),N_H(L),C_K(L))$ and $(C_G(L),C_H(L),C_K(L))$ are subfolders. As $\langle C_K(L) \rangle \le C_G(L)$
and $N_G(L) = N_H(L) \langle C_K(L) \rangle$ we have that $|N_G(L):C_G(L)|$ divides $|H|$.
By coprime action we have 
 $|N_{\overline{G}}(\overline{L})|_{2'} = |N_G(L)|_{2'}$, 
 $|C_{\overline{G}}(\overline{L})|_{2'} = |C_G(L)|_{2'}$ and $|H|_{2'} = |\overline{H}|_{2'}$, which implies the lemma.  
\end{bew}

\section{Passive simple groups: the classification} 
In this section $(G,H,K)$ is a BX2P-folder,
$\overline{G}=G/O_2(G)$, $F^\ast(\overline{G}) \cong S$ with $S$ some finite simple nonabelian group, 
$S \le T \le \Aut(S)$ with $\overline{G} \cong T$ and  $G_0$ the preimage of $F^\ast(\overline{G})$. 
Remember, that as a  starting point by \refNloopEmbedding, $\overline{H} \cap \overline{G_0}$ contains nontrivial elements of odd order. 

The goal of this section is : 

\begin{thm}
\label{MainClassificationTheorem}
Let $(G,H,K)$ be a BX2P-folder and $\overline{G}=G/O_2(G)$.
Suppose $F^\ast(\overline{G})$ is a finite nonabelian simple group. 
If $\overline{G} \ne \overline{H}$, then: 
\begin{itemize}
\item There exists an integer $q$ with $q=9$ or $q$ a Fermat prime, $q \ge 5$,
\item $\overline{G} \cong \PGL_2(q)$ or $\overline{G} \cong \PGaL_2(q)$ (for $q=9$ only),
\item $|\overline{G}:\overline{H}|=q+1$ and
\item $\overline{K}$ consists of 1 and all involutions in $\PGL_2(q)$ outside $\PSL_2(q)$. 
\end{itemize} 
In particular the nonpassive simple groups are among the groups $\PSL_2(q)$ for $q=9$ or a Fermat prime $q\ge 5$.
\end{thm}

\subsection{The groups $\PSL_2(q)$}
\begin{lemma}
\label{PSL2}
Let $S \cong \PSL_2(r), r>3$. If $\overline{G} \ne \overline{H}$, then 
$\overline{G} \cong \PGL_2(q)$ for $q=9$ or $q$ a Fermat prime
or $\overline{G} \cong \PGaL_2(9)$, $|G:O_2(G)H|=q+1$ and $\overline{K}$ consists of 1 and
all involutions of $\PGL_2(q)$ outside $\PSL_2(q)$. 
\end{lemma}

\begin{bew}
Let $r = p^e$ with $p$ a prime. Suppose first $p$ is odd.
By \refNloopEmbedding we get some element $x$ of odd order into $\overline{H} \cap \overline{G}^{(\infty)}$.
Suppose first, that $x$ is a $p'$-element, so $x$ is contained in some torus of size $\frac{r-1}{2}$ or $\frac{r+1}{2}$. 
Remember, that $\Aut(\PSL_2(r))$ has the following types of involutions: 
those in $\PSL_2(r)$, those in $\PGL_2(r)$ outside $\PSL_2(r)$ and possibly field automorphisms of order 2. 
The first two types of involutions invert both tori, so invert some conjugate of $x$.
So by \refnoHoverlineinvert $K$ cannot contain involutions of $\PGL_2(r)$, so consists of $1$ and field automorphisms only. 
Since field automorphism act nontrivially on a Sylow-$p$-subgroup, in this case $\overline{H}$ is a $p'$-group. 
We can now estimate the size of $\overline{K}$ and $|\overline{G}:\overline{H}|$: 
Let $r=s^2$. Then $|K| \le 1 + s (s^2+1)$. On the other hand $|G:H| \ge \frac{1}{2} s^2 (s^2-1)$. This gives a contradiction
since $s \ge 3$. So $\overline{K}^\sharp$ cannot consist of field automorphisms only or contain $p'$-elements of odd order. 
So $x$ is a $p$-element. Since $C_G(x)$ is soluble, but contains a Sylow-$p$-subgroup $P$ of $G$ we may assume by \refsolublegroups,
that $P \le H$. The Borel subgroup $N_{\overline{G}}(\overline{P})$ is then covered by $N_{G}(P)$
and $|G:O_2(G)H|=r+1$.
As $\overline{H}$ does not contain $p'$-elements of odd order, $r-1$ is a 2-power.
Notice, that in the case of $\overline{G}=\Aut(\Alt_6)= \PGaL_2(9)$ we still get $|\overline{G}:\overline{H}|=r+1$,
since the normalizer of a Sylow-3-subgroup of $\overline{G}$ has index 10 and is contained in $\overline{H}$.
Furthermore, in this case $\overline{K} \subseteq \PGL_2(q)$ since the other involutions are in $\Sigma_6$
and invert elements of order 3, which now cannot be in $\overline{K}$ by \refnoHoverlineinvert. 
So let $r$ be even, so $r\ge 4$. There are only two types of involutions, field automorphisms and inner automorphisms.
Inner automorphisms invert conjugates of all elements of odd order, so cannot be in $\overline{K}$. 
Field automorphisms act on a torus of size $r-1$ inside some invariant Borel subgroup, 
so $\overline{H}$ has to be the normalizer of a torus of size $r+1$. 
Calculation as in the case $r$ odd gives: $|\overline{K}| \le 1 + s (s^2+1)$ and 
$|\overline{G}:\overline{H}| \ge \frac{1}{2}s^2 (s^2-1)$, a contradiction for $s \ge 4$. 
The case $s=2$ was already handled as $\Alt_5 \cong \PSL_2(4)\cong \PSL_2(5)$.  
\end{bew}

\subsection{The alternating groups}
\begin{lemma}
\label{Alternating}
Let $S \cong \Alt_n$ for $n \ge 7$. Then $\overline{G}=\overline{H}$. 
\end{lemma}

\begin{bew}
Remember, that $\ell_T$ has property $FS_p$ for $p \ge 5$.  Furthermore \ref{alternating_centralizers} turns
out to be useful.

Let $n=7$. If $O_2(\overline{H})=1$, by \refHeissPrime, $7 \in \pi(H)$. By \ref{H_covers_automizer}, $7 \in \pi(H)$ implies $3 \in \pi(H)$,
which implies $3$-cycles in $\overline{H}$ by \ref{alternating_centralizers} and therefore a full Sylow-3-subgroup. No proper maximal subgroup
of $\overline{G}$ exists with this property by \cite{ATLAS}, so $\overline{H}=\overline{G}$ in this case.   
So $O_2(\overline{H}) \ne 1$. If $\overline{H}$ contains elements of order 3, then a full Sylow-3-subgroup by \ref{alternating_centralizers}.
A full Sylow-3-subgroup of $\Alt_7$ does not normalize any 2-subgroup of $\Sigma_7$ by \cite{ATLAS}, a contradiction. 
Again elements of order 7 in $\overline{H}$ imply elements of order 3 in $\overline{H}$  by \ref{H_covers_automizer}. 
So by \refTNloopEnvelope, $\overline{H}$ is a $\{2,5\}$-group with index at least $2 \cdot 3^2 \cdot 7= 126$.  
But elements of order 5 in $\Sigma_7$ are inverted by involutions, which are products of two or 3 commuting cycles, 
so we get $|\overline{K}|\le 1 + 21 < 126 \le |\overline{G}:\overline{H}|$ by \refnoHoverlineinvert, a contradiction.  

Let $n=8$. By \refHeissPrime we get $7 \in \pi(\overline{H})$ or $O_2(\overline{H})\ne 1$. 
If $7\in \pi(H)$, then $3 \in \pi(H)$ and $\overline{H}$ contains elements of order 3, which are the product of two 3-cycles.
By \ref{alternating_centralizers} then $\overline{H}$ contains $3$-cycles, so a Sylow-3-subgroup. By \cite{ATLAS} this implies, 
that $\overline{H}$ contains a subgroup isomorphic to $\Alt_7$, in which case $|\overline{G}:\overline{H}|$ is a 2-power, 
which implies, that the loop is soluble, so $\overline{H}=\overline{G}$.
So $O_2(\overline{H}) \ne 1$ and $7 \not\in \pi(H)$. If $H$ contains some element of order 3, which 
is a product of two 3-cycles, $H$ contains a Sylow-3-subgroup of $\overline{G}$. So $\overline{H}$ contains 3-cycles.
If $\overline{H}$ contains $3$-cycles, the centralizer of a 3-cycle contains $\Alt_5$, so $\overline{H}$
contains elements of order 5. Conversely if $\overline{H}$ contains elements of order 5, its centralizer contains  a normal 3-group
generated by a  3-cycle, so $\overline{H}$ contains 3-cycles. So $\overline{H}$ contains a subgroup
of order 15. No such subgroup $\overline{H}$ with $O_2(\overline{H}) \ne 1$ of $\Sigma_8=\Aut(\Alt_8)$ exists. 

Finally for $n\ge 9$ let $X \le H$ be a $p$-group for some odd prime $p$. By \ref{alternating_centralizers} we have $p$-cycles in $\overline{H}$. 
Furthermore for $p>3$ we have a full Sylow-$p$-subgroup in $\overline{H}$ by \ref{FS_p_alternating}, while for $p=3$ we have
$\overline{H}=\overline{G}$ by \ref{terminal_three_cycles}. If $n-p \ge 6$, the centralizer of a $p$-cycle
has a component of degree at least 6, so this component ends up in $\overline{H}$, and contains $3$-cycles.
If $n-p=5$, the index $|C_{\overline{G}}(\overline{x}):C_{\overline{H}}(\overline{x})|$ may be 6, but $C_{\overline{H}}(\overline{x})$ 
contains elements of order 5. Now $5$-cycles in $\overline{H}$ imply $3$-cycles in $\overline{H}$ for $n \ge 11$,  but also for $n=9$. 
In case $n=10$, $\overline{H}$ contains a Sylow-5-subgroup, so $O_2(\overline{H} )=1$ and by \refHeissPrime $3 \in \pi(H)$. 
If $n-p=3$ or $n-p=4$, we get 3-cycles into $\overline{H}$, since the centralizer is soluble. 

This leaves $n = p$, $n=p+1$ or $n=p+2$ for a prime $p$. Now if $p$ is not a Fermat prime, we get another odd prime $r$
dividing $p-1$ by \ref{H_covers_automizer}. 
Therefore $p$ is a Fermat prime $p\ge 17$ and $p$ is the unique odd prime dividing $|H|$. 
If $n=p$, then $\overline{H}$ is the normalizer of a Sylow-$p$-subgroup, which is a maximal subgroup of $\overline{G}$.
By \refNloopEmbedding we need a $\PGL_2(p)$ in $\overline{G}$ for a nonsoluble loop. As the permutation degree
of $\PGL_2(p)$ is $p+1$, we get $\overline{H}=\overline{G}$. 
If $n>p$, $\overline{H}$ cannot act transitively, so $\overline{H}$ is contained in the stabilizer of the
orbit decomposition. This stabilizer leads to a subloop by \refHsuper. By induction however $\overline{H}$
contains elements of order 3.
Since our arguments in this last case are based on the $N$-loop theorem of Aschbacher
in \cite{A} and it's counterpart in \cite{AKP}, it should be mentioned, that the direct approach of counting the involutions in $\overline{G}$ and 
comparing their number with the index of $\overline{H}$ was Aschbachers original argument. 
\end{bew}
\subsection{The sporadic groups}
\begin{lemma}
The sporadic simple groups are passive.
\end{lemma}

\begin{bew}
By \refHisTGiffGisTG, $\overline{H}$ is not a 2-group. 
By \ref{sporadiccyclicnormalizers} and \ref{H_covers_automizer} we may assume, that $\overline{H}$ contains an element of order $p$
with $p$ a Fermat prime, so $p=3,5$ or $17$. 

If $\overline{G}$ has only one class of involutions, the embedding of a 2M-loop by \refNloopEmbedding
shows, that involutions from this class invert some element of odd order in $\overline{H}$, a contradiction
to \refnoHoverlineinvert. Therefore $\overline{H}=\overline{G}$. For this reason $M_{11}$,$J_1$, $M_{23}$, $Ly$ and $Th$ 
are passive. 

Remember, that we have $FS_p$-property except $(p,S)$ is one of $(3,M_{23})$, $(3,HS)$  or $(5,Suz)$ by \ref{FS_p_sporadic}.    
We use the character tables in \cite{ATLAS} as provided in GAP for calculation of structure constants.
Specifically we calculated, which classes of involutions invert elements from classes of Fermat prime order.
For the structure of centralizers of elements we use without further reference the informations from \cite{ATLAS} 
in the list of maximal subgroups as well as the size of the centralizers from the character tables.

In case of $M_{12}$, structure constant calculations show, that $\overline{H}$ does not contain
elements of classes $3B$ or $5A$, as these classes are inverted by all classes of involutions. 
By $FS_3$ and $FS_5$-property, $\overline{H}$ does not contain elements of order 3 or 5, a contradiction to \refNloopEmbedding.

In case of $M_{22}$, elements from all conjugacy classes of involutions invert class $3A$, so $\overline{H}$
does not contain elements of order 3. From the list of maximal subgroups we conclude,
that $\overline{H}$ is contained in a maximal subgroup $\overline{M}$ of type $2^5: \Sigma_5$. All other classes
of maximal subgroups imply elements of order 3 in $\overline{H}$ by \refEnvelopeGroup. Furthermore
$\overline{K}$ consists of class $1A$ and $2B$, as elements from $2A$ and $2C$ invert elements of class $5A$. 
Now $|\overline{G}:\overline{H}| \ge 2 \cdot 3^2 \cdot 7 \cdot 11 = 1386  > |\overline{K}|=1 + 330 $, a contradiction.

In case of $J_2$, we get the following implications for containement in $\overline{H}$:
We have $FS_3$ and $FS_5$-property. Furthermore $H$-intersection with $3A$ implies intersection with $5AB$, 
while 5-elements in $H$ imply 3-elements in $H$ from the normalizer of a Sylow-5-subgroup. 
Among  maximal subgroups picked up by $\overline{H}$ are the normalizer of a $3A$-cyclic group and the normalizer of a Sylow-5-subgroup.
Therefore $\overline{H}=\overline{G}$.

In case of $HS$, elements from all classes of involutions invert elements from class $3A$, so $\overline{H}$
does not contain elements of order 3. (There is no class $3B$).
So $\overline{H}$ contains a Sylow-5-subgroup by $FS_5$-property. 
From a structure constant calculation we conclude, that $\overline{K}$ consists of classes $1A$ and $2C$.
A maximal subgroup containing a Sylow-5-subgroup is of type $U_3(5).2$ or the normalizer of the Sylow-5-subgroup (in $HS.2$).
By \refHsuper and \refEnvelopeGroup,  $\overline{H}$ is a $\{2,5\}$-group. 
Now $|\overline{G}:\overline{H}| \ge 2 \cdot 3^2 \cdot 7 \cdot 11 = 1386 > |\overline{K}| = 1+1100 $ gives a contradiction.

In case of $J_3$ we get from structure constant calculation, that $\overline{H}$
does not contain elements of order 3.
Elements of order 5 in $\overline{H}$ imply elements of order 3 in $\overline{H}$, since the centralizer of 5-elements is soluble
of size 30.
As $J_3.2$ does not involve a $\PGL_2(17)$, only a $\PSL_2(17)$, we conclude $\overline{H}=\overline{G}$ by
\refNloopEmbedding. 

In case of $M_{24}$ we get by structure constant calculations, that $\overline{H}$ does not contain any elements
of order 3 or 5, so $\overline{H}=\overline{G}$.

In case of $McL$, we can calculate that elements from class $2B$ (outer involutions) invert elements
from all classes of elements of order 3 and 5. So $\overline{K}$ does not contain outer involutions. 
But this gives a contradiction as in the case of $M_{11}$, $J_1$ etc. as above. 

In case of $He$, structure constant calculations show, that $\overline{H}$ does not contain elements of order 3.
If $\overline{H}$ contains elements of order 5, it contains a Sylow-5-subgroup. From the shape of the  normalizer of a Sylow-5-subgroup
we conclude, that then $\overline{H}$ contains elements of order 3. As $He$
does not contain elements of order $9$, $\PSL_2(17)$ is not involved in $\overline{G}$, so by \refNloopEmbedding 
$\overline{H}=\overline{G}$.

In case of $Ru$, structure constant calculations show, that $\overline{H}$ contains no elements of order 3
or 5, so $\overline{H}=\overline{G}$. 

In case of $Suz$, structure constant calculations show, that $\overline{H}$ does not contain elements from class $3C$ or $5B$.
By $FS_3$-property $\overline{H}$ does not contain elements of order 3. 
Furthermore $\overline{H}$ contains elements of order 3, if it contains elements of class $5A$, as the $5A$-centralizer may involve a $\PGL_2(9)$.
Therefore $\overline{H}=\overline{G}$. 

In case of $O'N$, structure constant calculations show, that $\overline{H}$ does not contain elements of order 3 or 5.

In case of $Co_3$, structure constant calculations show, that $\overline{H}$ does not contain elements of classes
$3B$, $3C$  or $5B$, so by $FS_p$-property no elements of order 3 or 5.

In case of $Co_2$, structure constant calculations show, that $\overline{H}$ does not contain elements of class $3B$,
so by $FS_3$-property no elements of order 3. 
Elements of class $5B$ in $\overline{H}$ imply elements of class $5A$ in $\overline{H}$, while the later imply elements of order 3 in $\overline{H}$. 

In case of $Fi_{22}$, elements of order 5 in $\overline{H}$ imply, that $\overline{H}$
contains a Sylow-5-subgroup of $\overline{G}$. From the list of maximal subgroups we conclude, that $\overline{H}$
is contained in a subgroup of type $O_8^+(2).\Sigma_3$ , ${}^ 2F_4(2)'$ ,$\Sigma_{10}$ or the corresponding maximal
subgroups in $Fi_{22}:2$. By \refEnvelopeGroup therefore $\overline{H}$ is among one of these groups.
Calculations of structure constants gives the bound $|\overline{K}|\le 65287$, if $\overline{H}$ contains elements of order 5.
Therefore $\overline{H}$ contains $\overline{M} \cong O_8^+(2).\Sigma_3$. This is a contradiction to $FS_3$-property, 
as the group in question does not contain a Sylow-3-subgroup of $G$. 
In particular we find in $\overline{H}$ a subgroup $\overline{P}$
of order $3^6$, such that $N_{\overline{H}}(\overline{P}) \not\le \overline{M}$. 
Therefore $\overline{H}$ does not contain elements of order 5 or $\overline{H}=\overline{G}$. 
Elements of class $3A$ in $\overline{H}$ imply elements of order 5 in $\overline{H}$, so $\overline{H}=\overline{G}$. 
By $FS_3$-property then $\overline{H}$ does not contain elements of order 3.

In case of $HN$, structure constant calculations show, that $\overline{H}$ does not contain
elements from classes $3A,5A$ or $5D$. By $FS_3$ and $FS_5$-property then $\overline{H}=\overline{G}$.

In case of $Fi_{23}$, structure constant calculations show, that $\overline{H}$ does not contain
elements of class $3A$, so by $FS_3$-property no elements of order 3. $FS_5$-property implies $K= 1A \cup 2A$, so
$|K| =31672$. But then $\overline{H}$ contains elements of order 3. So $\overline{H}$ does not contain elements of order 3 or 5. 
If $\overline{H}$ contains elements of order 17, then structure constant calculations show $|\overline{K}| \le 55614277$, but the index of a $\{2,17\}$-subgroup 
is at least $2 \cdot 3^{13} \cdot 5^2 \cdot 7\cdot 11\cdot 13 \cdot 23 = 1835304921450$. 

In case of $Co_1$, structure constant calculations show, that $\overline{H}$ does not contain elements of classes
$3B,3D$ or $5B$, so by $FS_p$-property no elements of order 3 or 5.  

In case of $J_4$, structure constant calculations show, that $\overline{H}$ does not contain elements of order 3 or 5. 

In case of $Fi_{24}'$, structure constant calculations show, that $\overline{H}$ does not contain elements
of class $3A$, so no elements of order 3 by $FS_3$-property.  
Structure constant calculations also show, if $\overline{H}$ contains elements of order 5, $\overline{K}$ consists
of $1A$ and $2C$ only, so $|\overline{K}|=1+306936$. In that case $\overline{H}$ contains a subgroup isomorphic to $Fi_{23}$,
so contains elements of order 3. 
Now $\overline{H}$ is a $\{ 2,17\}$-group with $|\overline{K}| \le 4860791965$. (Bound obtained from structure constant calculations.)
But $|\overline{G}:\overline{H}| \ge 2 \cdot 3^{16} \cdot 5^2 \cdot 7^3 \cdot 11 \cdot 13 \cdot 23 \cdot 29$ = 70415143921272150,
a contradiction. 

In case of $B$, structure constant calculations show, that $\overline{H}$ does not contain elements of class $3A$,
so by $FS_3$-property no elements of order 3. 
Elements of classes $5A$ or $5B$ in $\overline{H}$ imply elements of order 3 in $\overline{H}$, so by $FS_5$-property
$\overline{H}$ does not contain elements of order 5. 
Again $\overline{H}$ is a $\{2,17\}$-subgroup. Structure constant calculations show $|\overline{K}|\le 11721020628376$,
but $|\overline{G}:\overline{H}| \ge \frac{|B|}{2^{40} \cdot 17} = 222279514364689031250$. 

Finally in case $M$, structure constant calculations show, that $\overline{H}$ does not contain elements
of classes $3A, 3C$ or $5A$, so no elements of order 3 or 5 by $FS_p$-property.
Elements of class $17A$ in $\overline{H}$ imply elements of order 3 in $\overline{H}$. 
\end{bew}
\subsection{Groups of Lie type in odd characteristic}
Let $S$ be a simple group of Lie type in characteristic $p>2$ and 
$q=p^f$ be the field parameter of $S$, $q$ odd in this paragraph. 
(If $S$ is defined relative to a field extension, $q$ is the size of the smaller field.)
Before using \ref{OddChar} we should handle the cases $S \cong \PSU_3(q)$ and $S \cong {}^2G_2(q)$: 

\begin{lemma}
\label{PSU3_odd}
Let $S \cong \PSU_3(q)$ for $q$ odd. Then $\overline{G}=\overline{H}$. 
\end{lemma}

\begin{bew}
Let $d:=(q+1,3)$ and $s \in \pi(H)$ with $s$ odd. By \refNloopEmbedding such an $s$ exists.
Let $\overline{x} \in \overline{H}$ be of order $s$. If $s$ divides $\frac{q^2-q+1}{d}$, 
then we get $3 \in \pi(H)$ by \ref{H_covers_automizer} and $3$ divides $(q-1)q(q+1)$. So we may assume
from the start, that $s$ divides $q-1$, $q$ or $q+1$. 

Notice, that $\Aut(S)$ has only two classes of involutions: inner involutions of $S$ with centralizer $\ZZ_\frac{q+1}{d} \cdot \SL_2(q):2$
and outer involutions with centralizer $O_3(q) \times \ZZ_2 \cong \PGL_2(q) \times \ZZ_2$. (See 4.5.1 of \cite{GLS3} for details.)
Remember that by \ref{Out_restriction} $T/S$ is a subgroup of $\Sigma_3$ not of order 3.
Let $i \in \Aut(S)$ be an outer involution with $C=C_S(i) \cong \PGL_2(q)$. As $S$ has only one class of involutions, elements
of order $p$, $q+1$ and $q-1$ are inverted by inner involutions of $S$. As $\Aut(S)$ has only one class of outer involutions,
so all involutions of $C_T(i)-C_S(i)$ are conjugate to $i$, also outer involutions invert elements of order $p$, $q+1$ and $q-1$. 

By \ref{FS_characteristic_p}, $\overline{H}$ does not contain any elements of order $p$, as otherwise $\overline{H}$
would contain a Sylow-$p$-subgroup of $\overline{G}$ and some element of $\overline{K}$ would invert some element of $\overline{H}$ of odd order $p$.

If $s$ divides $q-1$, then $C_{\overline{G}}(\overline{x})$ is soluble, so $\ell_T$ has property $FS_s$,
and as in case $s=p$ we get a contradiction. 

This leaves the case, that $s$ divides $q+1$. Then either $C_{\overline{G}}(\overline{x})$ is soluble or contains a unique
$\SL_2(q)$-component. The later case would imply elements of order $p$ in $\overline{H}$. So $C_{\overline{G}}(\overline{x})$
is soluble. If the Sylow-$s$-subgroup is abelian, $\overline{H}$ contains a Sylow-$s$-subgroup, so elements, which are inverted by some involution of $K$, 
a contradiction to \refnoHoverlineinvert. 

The only remaining possibility is $s=3$ and $3|q+1$. But then $N_{\overline{G}}(O_3(C_{\overline{G}}(\overline{x}))) \le \overline{H}$
and $\overline{H}$ contains a Sylow-$3$-subgroup of $\overline{G_0}$, which is again a contradiction, as every involution of $\Aut(S)$ inverts some element of order $3$ in $S$. 
\end{bew}

\begin{lemma}
\label{Ree2G2}
Let $S \cong {}^2G_2(q)$. Then $\overline{G}=\overline{H}$.  
\end{lemma}

\begin{bew}
Remember that $\Aut(S)$ does not contain outer involutions and only one class of inner involutions.
Then \refNloopEmbedding gives a contradiction, as $\overline{H}$ contains involutions, which invert nontrivial elements of odd order in $\overline{H}$. 
\end{bew}

Let $\overline{x} \in \overline{H}\cap \overline{G_0}$ be some element of odd prime order $r$. 
By \ref{OddChar} we may assume that $r \ne p$ and  $C_{\overline{G}}(\overline{x})$ is a soluble $p'$-group
as otherwise $\overline{H}=\overline{G}$. 

\begin{lemma}
\label{PSL3_odd}
Let $S \cong \PSL_3(q)$ for $q$ odd. Then $\overline{H}=\overline{G}$.
\end{lemma}

\begin{bew}
We handle the two cases separately: $q$ a square and $q$ not a square.
If $q$ is not a square, then by Theorem 4.5.1 of \cite{GLS3} $\Aut(S)$ has only two classes of involutions,
inner involutions and graph automorphisms of order 2, which centralize a $\PGL_2(q) \cong O_3(q)$ in $S$. 
We see, that inner involutions invert elements of order $p$, $q-1$ and $q+1$. But in the direct product
$\PGL_2(q) \times \ZZ_2$, which is isomorphic to the centralizer of graph automorphism $j$, the products
$ij$ with $i$ an inner involution of $S$, $i \in \PGL_2(q)$, all are outer involutions, but invert the same
elements as $i$. Therefore any involution of $\Aut(S) $ inverts some elements of order $p$,$q-1$ and $q+1$.
By \refnoHoverlineinvert and $FS_r$-property, $\overline{H}$ is therefore not divisible by $p,q-1$ or $q+1$.
By \refNloopEmbedding $\overline{H} \cap \overline{G_0}$ contains elements of odd order $r$, so $r$ divides
$q^2+q+1$. Now the centralizer of an element $\overline{x} \in \overline{H}$ of order $r$ is cyclic
and by \ref{H_covers_automizer} we get  $3 \in \pi(H)$, a contradiction as $3$ divides $(q-1)q(q+1)$.
If $q=q_0^2$ is a square, we have four classes of involutions: in addition to inner involutions and graph automorphisms
we get field automorphisms and graph-field automorphisms into $\Aut(S)$. Field automorphisms
centralize a $\PSL_3(q_0)$, while graph-field automorphisms centralize a $\PSU_3(q_0)$. Notice, that $3$ divides
$(q_0-1)q_0(q_0+1)$. Let $r \in \pi(H)$, $r$ odd. We show, that $r \nmid q_0^2-1$: 
As in case $q$ not a square, graph and inner involutions invert elements of order $q-1 = q_0^2-1$. 
But also field and graph field involutions invert elements of order $q_0-1$ and $q_0+1$ with the same argument
as in the case $q$ not a square. Notice, that the involutions $i$ and $ij$ with $[i,j]=1$, $i$ an inner involution
and $j$ an outer involution, are in the same  $S$-coset. By Theorem 4.9.1 of \cite{GLS3} these involutions are conjugate.
By \refnoHoverlineinvert and $FS_r$-property, $\overline{K}$ does not contain involutions of $\overline{G}$,
so $\overline{H}=\overline{G}$, if $r \mid q_0^2-1=q-1$. 
If $r \mid q^2+q+1 = (q_0^2+q_0+1)(q_0^2-q_0+1)$, then by \ref{H_covers_automizer} we have $3 \in \pi(H)$,
but $3 \mid (q_0-1) q_0 (q_0+1) \mid (q-1) q $, so $\overline{H}=\overline{G}$ in this case.
So remains $r \mid q_0^2+1 = q+1$. Then $K$ consists of 1 and field and/or graph-field involutions. 
The index $|\overline{G}:\overline{H}|$ is divisible by 2, $q^3 =q_0^6$ and $q^2+q+1=q_0^4+q_0^2+1$. 
On the other hand $|\overline{K}| \le 1 + q_0^3(q_0^2+1)(q_0^3-1) + q_0^3(q_0^2+1)(q_0^3+1)= 1+ 2 q_0^6 (q_0^2+1)$,
which gives the contradiction $|\overline{K}|<|\overline{G}:\overline{H}|$.  
\end{bew}

\begin{lemma}
Let $S \cong \PSL_n(q)$ or $\PSU_n(q)$ with $n \ge 4$, $q$ odd.
Then $\overline{H}=\overline{G}$.
\end{lemma}

\begin{bew}
By \ref{OddChar}, we may assume, that $\overline{x}$ is not in the big connected component of $\Gamma_{\cal O}$.
We use \cite{CommGraph}, Thm 4 for a list of small connected components.

In case $n=4$ and $d_q(r) \in \{3,6\}$  we use \ref{H_covers_automizer} to get $3 \in \pi(H)$,
but elements of order 3 are in the big connected component.

In case $n=4$ and $d_q(r)=4$ we determine the the structure of maximal subgroups $\overline{M}$ of $\overline{G}$, using \cite{KL}. 
By \refEnvelopeGroup we conclude that either $3 \in \pi(H)$ or $p \in \pi(H)$ with elements of order 3 in the big connected component.

In case $n \ge 5$, $n-1$ a prime we use \ref{H_covers_automizer} to get $n-1 \in \pi(H)$. 
But elements of order $n-1$ are in the big connected component, as $d_q(n-1) | n-2$. 

In case $n \ge 5$ and $n$ a prime we use \ref{H_covers_automizer}, so $n \in \pi(H)$. 
We have $d_q(n)\mid n-1$. If elements of order $n$ are in the big connected component, then $3 \in \pi(H)$, 
so elements of order $n$ are in a another small connected component and therefore $n=5$.

Together with the remaining small connected components for $n=5$ this gives either $\pi(H) \subseteq \{2,5,11 \}$ in $\PSL_5(3)$ or  
$\pi(H) \subseteq \{2,5,61 \}$ in $\PSU_5(3)$.
In $\PSL_5(3)$, if $11 \in \pi(H)$, then $\overline{H}$ contains a torus normalizer $\ZZ_{11}^2:\ZZ_5$, which is a maximal subgroup 
of $\overline{G_0}$.
In $\PSU_5(3)$, if $61 \in \pi(H)$, then $\overline{H}$ contains a torus normalizer $\ZZ_{61}:\ZZ_5$, which is a maximal subgroup 
of $\overline{G_0}$.
In both cases we get a contradiction, as elements of order $5$ should be inverted in $\overline{H} \cap \overline{G_0}$ by \refNloopEmbedding.

So  $\overline{H}$ is a $\{2,5\}$-group. We use the list of maximal subgroups in \cite{KL} to determine
the possible $\overline{M}$ containing $\overline{H}$. In almost all cases then \refEnvelopeGroup produces
additional primes into $\pi(H)$. 
The only remaining case is a subgroup of type $\ZZ_4^5 : \Sigma_5$ in $\PSU_5(3)$. 
In this case elements of order 5 are inverted by involutions from all but one conjugacy class,
a class of length 4941. Since $|\overline{G}:\overline{H}| \ge 2 \cdot 3^{10} \cdot 7 \cdot 61$,
$\overline{H}=\overline{G}$ in this case. 
\end{bew}

\begin{lemma}
Let $S \cong \PSp_{2n}(q)$ for $n \ge 2$ or $\POmega_{2n+1}(q)$ for $n \ge 3$ with $q$ odd.
Then $\overline{H}=\overline{G}$.
\end{lemma}

\begin{bew}
By \ref{OddChar}, we may assume, that $\overline{x}$ is not in the big connected component of $\Gamma_{\cal O}$.
We use \cite{CommGraph},Thm 4 for a list of small connected components.
If $q-1$ or $q+1$ is a 2-power and $n$ is an odd prime, we use \ref{H_covers_automizer} to get $n \in \pi(H)$. 
As $d_q(n) | n-1$, we have elements of order $n$ in the big connected component, so $\overline{H}=\overline{G}$. 
In the remaining cases $n$ is a 2-power and $r \mid q^n+1$. We determine the structure of possible maximal
subgroups $\overline{M}$, which contain $\overline{H}$, using the list in \cite{KL}. 

In the symplectic case we get candidates:
$\overline{M_1}$ of type $\Sp_n(q^2).2$ in class ${\cal C}_3$, 
$\overline{M_2}$ of type $\GU_n(q)$ in class ${\cal C}_3$, 
$\overline{M_3}$ of type $2^{1+2m}.O_{2m}^-(2)$ in class ${\cal C}_6$,
or $\overline{M_4}$ in class ${\cal S}$ with $F^\ast(\overline{M_4})$ a simple group.
By \refHsuper and \refEnvelopeGroup, $\overline{M_1}$ and $\overline{M_2}$ imply $p \in \overline{H}$. 
In case $\overline{M_3}$ we have $2n=2^m$ and the largest prime dividing $|O^-_{2m}(2)|$ is 
bounded by $2^m+1=2n+1$. On the other hand, for each odd prime $r$ dividing $|H|$ we have $d_q(r) =2n$,
so $r \ge 2n+1$. Therefore $\overline{M_3}$ contains the torus iff $\frac{q^n+1}{2}=2n+1$, which holds
only for $q=3$, $n=2$.  
In case $M_4$, if $F^\ast(\overline{M})$ is passive and not a Suzuki group, then $3 \in \pi(H)$ and elements of order 3
are in the big connected component. If $\overline{M_4}$ is a Suzuki group, then $5 \in \pi(H)$, but by Landazuri-Seitz, $n \ge 4$, so elements of order 5 are in the big connected component.
Remains $F^\ast(\overline{M}) =\PSL_2(q_1)$ for $q_1=9$ or $q_1$ a Fermat prime. As then $\overline{H} \cap \overline{M}_4$
contains a Borel subgroup, we have $\frac{q^n+1}{2}=q_1$.  
Suppose $q_1 =2^e+1$ with $e$ even. Then $\frac{q^n+1}{2} = 2^e+1$ gives
$q^n-1 = 2^{e+1}$, so $q=3,n=2,e=2$, which is again the special case of $\Sp_4(3)$. 

In case $\PSp_4(3)\cong \PSU_4(2)$ we still have to exclude the case, that $\overline{H}$ is a $\{2,5\}$-group. We use information from \cite{ATLAS}.
The only possible subgroup $\overline{M}$ is of shape $2^4:\Sigma_5$ and of index 27. This gives a candidate
for $\overline{H}$ of index $6 \cdot 27$.
Involutions not inverting elements of order 5 are in class 2A and 2C, so $|\overline{K}| \le 1+45+36 = 82$.
As $|\overline{G}:\overline{H}| \ge 2 \cdot 3^4 = 2 \cdot 81$ we get a contradiction. 

In the orthogonal case we get a subgroup of type $O_1(q) \perp O_{2n}^-(q)$ and maybe maximal subgroups in class {\cal S}. 
From \refEnvelopeGroup we conclude, that $\overline{H}$ contains other elements of odd order, so $\overline{H}=\overline{G}$ or a maximal subgroup is of type $\PGL_2(q_1)$ for $q_1\ge 257$ a Fermat prime. 
(Since $n \ge 3$, $\frac{q^n+1}{2}$ is at least $\frac{3^4+1} {2}=41$.) 
Now $q_1=2^e+1$ for $e$ even, so if $\frac{q^n+1}{2} = 2^e+1$, then $q^n-1= 2^e$, which happens only for $q=3$, $n=2$, 
a case which is excluded by $n \ge 3$.
\end{bew}

\begin{lemma}
Let $S \cong \POmega^+_{2n}(q)$ or $\POmega^-_{2n}(q)$ for $n \ge 4$, $q$ odd.
Then $\overline{H}=\overline{G}$.
\end{lemma}

\begin{bew}
By \ref{OddChar}, we may assume, that $\overline{x}$ is not in the big connected component of $\Gamma_{\cal O}$.
We use \cite{CommGraph},Thm 4 for a list of small connected components.
If $n$ is an odd prime, by \ref{H_covers_automizer} $n \in \pi(H)$.
From the list of small connected components, we conclude, that $n$ is in the big connected component.
If $n-1$ is an odd prime, $n-1 \in \pi(H)$ by \ref{H_covers_automizer}. Again $n-1$ is in the big connected component.
In the remaining cases we use the list of maximal subgroups in \cite{KL} for possible maximal subgroup $\overline{M}$
containing $\overline{H}$. The following observations eliminate some maximal subgroups:

If $n-1$ is a 2-power, as $d_q(r) = 2n-2$, either $r=2n-1$ or $r \ge 2(2n-2)+1 = 4n-1$, so $r$ in the big connected component. 
If $r=2n-1$, then $r$ is a Fermat prime. If $\frac{q^{n-1}+1}{2}=2n-1$, then $q=3$ and $n=3$ contrary to $n \ge 4$.
In case $n$ a 2-power, a similiar arguments works.

Therefore, if $\overline{M}$ is not in class ${\cal S}$, we have:
In $\Omega^+_{2n}(q)$, if $(n-1)(q+1)$ is a 2-power and $d_q(r)=2n-2$, then $\overline{M}$ is of type $O_2^-(q) \perp O_{2n-2}^-(q)$ in class ${\cal C}_1$ or 
a subgroup of type $O_n(q^2)$ in class ${\cal C}_3$.
In $\Omega^-_{2n}(q)$, if $(n-1)(q-1)$ is a 2-power and $d_q(r)=2n-2$, then $\overline{M}$ is of type $O_2^+(q) \perp O_{2n-2}^-(q)$ in class ${\cal C}_1$,
a parabolic subgroup of type $P_1$ also from class ${\cal C}_1$, a subgroup of type $\GU_n(q)$ in class ${\cal C}_3$ or  a subgroup of type $O_n(q^2)$ also from ${\cal C}_3$.
In $\Omega^-_{2n}(q)$, if $n$ is a 2-power and $d_q(r)=2n$, then $\overline{M}$ is of type $O^-_n(q^2)$ in class ${\cal C}_3$.
By \refHsuper and \refEnvelopeGroup then either $3 \in \pi(H)$ or $5 \in \pi(H)$ with elements of order 3
and 5 in the big connected component or $\overline{M}$ is of type $\PSL_2(q_1)$ for a Fermat prime $q_1 \ge \frac{3^4+1}{2}=41$, so $q_1 \ge 257$. We get a contradiction as a faithful representation of $\PSL_2(q_1)$ has degree
at least $\frac{q_1-1}{2}$, so $\frac{q_1-1}{2} \le 2n$, but $q_1 \ge \frac{q^{2n-2}+1}{2}$.
\end{bew}

\begin{lemma}
Let $S$ isomorphic to one of  
$G_2(q),$ ${}^3D_4(q)$, $F_4(q),$ $E_6(q),$ ${}^2E_6(q),$ $E_7(q)$ or $E_8(q)$ for $q$ odd. Then $\overline{H}=\overline{G}$.
\end{lemma}

\begin{bew}
By \ref{OddChar}, we may assume, that $\overline{x}$ is not in the big connected component of $\Gamma_{\cal O}$.
We use \refsmallCC for a list of small connected components.
In every case, the exceptions come from self centralizing tori $T$ in $S$. If $|N_S(T):T|$ is not a 2-power,
we get $s \in \pi(H)$ for $s$ some odd prime divisor of $|N_S(T):T|$. Notice, that $s \in \{3,5,7\}$
and in all these cases the big connected component contains elements of order $s$, so $\overline{H}=\overline{G}$. 
So the following cases of $(S,d_q(r))$ remain:
$({}^3D_4(q),12)$, $(F_4(q),8)$, $(E_6(3),8)$, $(E_6(7),8)$, $({}^2E_6(3),8)$ or $({}^2E_6(7),8)$.
In these cases we have either $O_2(\overline{H})=1$ or $O_2(\overline{H}) \ne 1$.
If $O_2(\overline{H})=1$, the prime $p$ itself satisfies the prerequisites of \refHeissPrime, 
so by \refHeissPrime and \ref{OddChar}:  $\overline{G}=\overline{H}$.
Else $O_2(\overline{H}) \ne 1$, so $\overline{H}$ is a 2-local subgroup of $\overline{G}$. 
Now \cite{CLSS} and \cite{LSS} give a list of maximal subgroups, which contain all maximal local subgroups except
centralizers of outer automorphisms. The structure of centralizers of outer involutions
in these cases is described in \cite{GLS3}. So for any $\overline{H}$ we know the structure of at least one maximal subgroup $\overline{M}$
containing $\overline{H}$. 
By \refHsuper, we can use \refEnvelopeGroup on maximal subgroups $\overline{M}$, which contain $\overline{H}$. 
As a result we get elements of order $3$ into $\overline{H}$, so by \ref{OddChar} $\overline{H}=\overline{G}$. 

Notice, that in case of $\overline{G} \cong {}^3D_4(q)$ the list of maximal subgroups actually
produces the torus normalizer (of type $q^4-q^2+1:4$) itself as the unique maximal subgroup containing the torus.
In that case we conclude, that $O_2(\overline{H})=1$, since outer involutions are field automorphisms and act
on the torus nontrivially.
\end{bew}
\subsection{Groups of Lie type in even characteristic}
The main arguments used in this paragraph are \ref{FS_p_even_char}, \ref{EvenCharTerminals}, \ref{soluble_subloops_in_even_char}
and the results on the commuting graph $\Gamma_{\cal O}$ in \cite{CommGraph} and special centralizers from section 2.

We first handle groups of low rank, with subcases $q=2$, $q=4$ and $q>4$ and later
the generic case with subcases $q\ge 4$ and $q=2$. 
  
\subsubsection{Low rank}
Recall, that we handled already $\PSL_2(q)$ in \ref{PSL2} and that
$Suz(q)={}^2B_2(q)$ is passive due to $2M$-loop embedding, \refNloopEmbedding.

In this section we handle the cases $S$ of type $\PSL_3(q)$, $\PSL_4(q)$, $\PSU_3(q)$, $\PSU_4(q)$, $\Sp_4(q)$, $G_2(q)$, ${}^3D_4(q)$, 
and ${}^2F_4(q)$. 
Let $S \le T \le \Aut(S)$ and $G/O_2(G) \cong T$. Remember the $FS_p$-property
from \ref{FS_p_even_char}. 
We will make the case division $q=2$, $q=4$ and $q>4$. 

\begin{lemma}
If $q=2$, then $\overline{H}=\overline{G}$. 
\end{lemma}
\begin{bew}
In case $q=2$ we can already exclude some groups for the following reasons: 
$\PSL_3(2)$ because of $2M$-Loop-embedding, \refNloopEmbedding,
$\PSL_4(2)$ because of the isomorphism with $\Alt_8$, 
$\PSU_3(2)$ because the group is soluble, 
$\PSU_4(2)$ because of the isomorphism with $\PSp_4(3)$, 
$\Sp_4(2)'$ because of the isomorphism with $\Alt_6 \cong \PSL_2(9)$,
$G_2(2)'$ because of the isomorphism with $U_3(3)$ and $2M$-Loop-embedding, \refNloopEmbedding.

The groups ${}^3D_4(2)$ and ${}^2F_4(2)'$ are ATLAS-groups, so we can use the
information from \cite{ATLAS}: By  $2M$-Loop-embedding and the list of maximal
subgroups we conclude, that if $S \cong {}^3D_4(2)$, then $\overline{H}=\overline{G}.$ 

For the Tits group ${}^2F_4(2)'$ we can establish the $FS_p$-property for all primes $p>2$
by \ref{FS_p_criterion}(1). 
Notice, that ${}^2F_4(2) = \Aut({}^2F_4(2)')$ is not generated by involutions. 
From the list of maximal subgroups in \cite{ATLAS} we conclude, that $\overline{M}$, 
a maximal subgroup of $\overline{G}$ which contains $\overline{H}$,
is isomorphic to $\PSL_3(3).2$. This implies $\overline{H}=\overline{M}$
and $O_2(\overline{H})=1$. By \refHeissPrime we get a contradiction, as 
the length of both classes $2A$ and $2B$ is divisible by $5$, so $\overline{H}$
has to contain a Sylow-5-subgroup of $\overline{G}$ too. 
\end{bew}

Remember, that for $q=4$ we have already the $FS_p$-property for all odd primes $p>3$. 

\begin{lemma}
If $q=4$, then $\overline{G}=\overline{H}$ or $S$ is of type $L_3(4)$.  
\end{lemma}

\begin{bew}
Let $S\cong \PSL_4(4)$. Notice, that if $5 \in \pi(H)$, then $\overline{H}$ 
contains a Sylow-5-subgroup. The normalizer of a 5-Sylow-subgroup contains
elements of order 3, while there exist elements of order $85=5 \cdot 17$ in
$G$, so then $|H|$ is divisible by $3 \cdot 5^2 \cdot 17$, and contains  
subgroups of type $5^2:3$ and $5 \times 17$. No such proper subgroup exists
by the list of maximal subgroups in \cite{KL}. If $5 \not\in \pi(H)$, 
then also $17 \not\in \pi(H)$ and elements of order 3 in $\overline{H}$
do not commute with elements of order 5 in $\overline{G}$. Else by the
structure of the centralizers of elements of order 3 we get $5 \in \pi(H)$ by 
\refEnvelopeGroup. No such elements of order 3 exist. Elements of order 7 in
$\overline{H}$ imply elements of order 3 in $\overline{H}$ both by the
centralizer of elements of order 7 and the normalizer of subgroups of order
7. We now get a contradiction to \refNloopEmbedding.

Let $S \cong \PSU_3(4)$. We use notation of p.30 of \cite{ATLAS}.
If $3 \in \pi(H)$, then $5 \in \pi(H)$ as the Centralizer of a $3A$-element is cyclic
of order 15. Furthermore the centralizer of a $5ABC-$ or $D$-element does not involve a $\PGaL_2(4)$, only
a $\PSL_2(4)$. So for $x \in \overline{H}$ of order 5 we have $O^2(C_{\overline{G}}(\overline{x})) \le \overline{H}$ by \ref{soluble_subloops_in_even_char}.
The normalizer of a Sylow-5-subgroup is a maximal subgroup of type $5^2: \Sigma_3$.
Therefore $\overline{H}=\overline{G}$, if $\overline{H}$ contains elements of order 3 or 5. 
Elements of order 13 in $\overline{H}$ imply elements of order 3 in $\overline{H}$. 

Let $S \cong \PSU_4(4)$. We calculated some centralizer data using MAGMA:
Elements of order 17 in $\overline{H}$ imply elements of order 3 in $\overline{H}$, as $\GU_4(4)$
contains a $\GL_2(16)$, so the centralizer of an element of order 17 is soluble and contains elements of order 3.
Elements of order 13 in $\overline{H}$ imply elements of order 3 in $\overline{H}$ by \ref{H_covers_automizer}
as already visible in $\PSU_3(4)$. 
Elements of order 3 in $\overline{H}$ imply elements of order 5 in $\overline{H}$: 
There are two classes of subgroups of order 3 with components $\PSL_2(16)$ and $\PSL_2(4)$ respectively in their centralizer.
In both cases $3 \in \pi(H)$ implies, that $5 \in \pi(H)$. 
If $5 \in \pi(H)$, by $FS_5$-property a Sylow-5-subgroup of $\overline{G}$ is already in $\overline{H}$.
There are elements of order 5 in that Sylow-subgroup, whose centralizer has shape $5 \times \PSU_3(4)$.
For these elements we can use \ref{soluble_subloops_in_even_char} to get the component $\PSU_3(4)$ into $\overline{H}$.
Calculation reveals, that there are so many elements of that type in a Sylow-5-subgroup, that $\overline{H}=\overline{G}$. 

Let $S \cong Sp_4(4)$. Elements of order 3 in $\overline{H}$ imply $5 \in \pi(H)$,
$5 \in \pi(H)$ or $17 \in \pi(H)$ implies $O_2(\overline{H})=1$. $O_2(\overline{H})=1$ implies
that $\overline{H}$ contains Sylow-5- and Sylow-17-subgroups of $\overline{G}$ in $\overline{H}$ by \refHeissPrime.
This implies $\overline{H}=\overline{G}$ by \cite{ATLAS}.

Let $S\cong G_2(4)$. We calculate centralizers with MAGMA, using the 12-dimensional representation of $G_2(4).2$
over $GF(2)$. In particular the subloops to centralizers of elements of order 5 are soluble,
as not sections $\PGaL_2(4)$ are involved, (only $\PSL_2(4)$.) So if $\overline{x} \in \overline{H}$ is of order 5,
then $O^2(C_{\overline{G}}(\overline{x})) \le \overline{H}$ by \ref{soluble_subloops_in_even_char}. 
Then $\overline{H}=\overline{G}$, as $\overline{H}$ does not only contains a $\PSU_3(4)$, but also
subgroups of type $5 \times A_5$ from both conjugacy classes, while $\PSU_3(4)$ contains only one such class.
Elements of order 3 in $\overline{H}$ imply elements of order 5 in $H$ by the centralizer structure, 
while elements of order 7 or 13 in $\overline{H}$ imply elements of order 3 in $\overline{H}$ by \ref{H_covers_automizer}.   

Let $S \cong {}^3D_4(4)$. If $3 \in \pi(H)$, then $5 \in \pi(H)$:  
Either $FS_3$-property fails on a subnormal $\PGaL_2(4)$ in a centralizer of an element of order 3
or $FS_3$-property holds. The first case implies $5 \in \pi(H)$, while
the second case implies elements of order 3 in $\overline{H}$ with centralizer shape
$(7 \times \SL_3(4)).3$, so again $5 \in \pi(H)$. 
By $2M$-loop embedding \refNloopEmbedding either $3 \in \pi(H)$ or $5 \in \pi(H)$, so $5 \in \pi(H)$
and $\overline{H}$ contains a Sylow-5-subgroup of $\overline{G}$. As there are centralizers of elements of order 5
of shape $5 \times \PSL_2(64)$, we have $3,7,13 \in \pi(H)$, so $\overline{H}$ contains Sylow-subgroups for the primes
$5,7$ and $13$. From the list of maximal subgroups of \cite{K3D4} we conclude, that $\overline{H}=\overline{G}$.  
\end{bew}

\begin{lemma}
Let $S \cong \PSL_3(4)$. Then $\overline{G} = \overline{H}$.
\end{lemma}

\begin{bew}
This group needs special treatment due to the exception of Zsygmondy's theorem and the fact, 
that $q-1 = (3,q-1)$.
We use Atlas-notation for the conjugacy classes, see \cite{ATLAS},p. 23.
Conjugacy classes of odd prime order are $3A, 5AB$ and $7AB$, of each odd prime order there
is a unique conjugacy class of groups of that order in $S \cong \PSL_3(4)$. 
$\Aut(\PSL_3(4))$ has involution conjugacy classes $2A,2B,2C$ and $2D$, which invert the following conjugacy classes of odd prime order:
$2A$ inverts elements from $3A$ and $5AB$. 
$2B$ inverts elements from $3A$ and $7AB$.
$2C$ inverts elements of all classes of 3-elements (including outer classes).
$2D$ inverts elements from $3A$, $5AB$ and $7AB$. 
Therefore $\overline{K}$ does not include elements of class $2D$.
As any involution inverts elements of class $3A$, $\overline{H}\cap \overline{G_0}$ does not contain elements of order 3.
So $\overline{H}\cap \overline{G_0}$ is a $\{2,5\}$-group by the $2M$-loop embedding, \refNloopEmbedding.
In particular maximal subgroups containing $\overline{H}$ have no $\Alt_6$ or $\PSL_3(2)$-components
and $\overline{K}$ does not contain involutions from $2A$ or $2D$.

Notice, that class $2B$ has length 280, while class $2C$ has length 120 or 360, depending on the presence
of diagonal automorphisms of order 3. Class $2B$ is a class of graph-field automorphisms, while class $2C$ is
a class of field automorphisms. 
We now check $\overline{G}$ for possible maximal subgroups $\overline{M}$ containing $\overline{H}$. 
Let $X_0:=\Aut(S) = L_3(4).D_{12}$. Calculations of maximal subgroups were done in MAGMA, using a 42-point representation
of $X_0$.

$X_0$ has 8 classes of maximal subgroups: 
$X_1 \cong \Sigma_3 \times \Sigma_5$, but if $\overline{M}$ is $X_1$, 
then $|\overline{G}:\overline{H}|\ge 6 \cdot |\overline{X_0} :\overline{X_1}|=2016$, while $|\overline{K}| \le 1+280+360=641.$
$X_2,X_3,X_4$ are soluble of sizes $2^2 \cdot 3^2 \cdot 7$,$2^5 \cdot 3^3$ resp. $2^8 \cdot 3^2$. 
$X_5 \cong \PSL_3(4).2^2$ is analyzed below, 
$X_6 \cong \PSL_3(4).6$ contains involutions of classes $2A$ and $2B$, but is not generated by involutions. 
See $X_{5,7}$ below for $\PSL_3(4).2$ with $2B$-outer involutions. 
$X_7 \cong \PSL_3(4).\Sigma_3$ with $2C$-outer involutions is analyzed below, but 
$X_8 \cong \PSL_3(4).\Sigma_3$ with $2D$-outer involutions is out.

$X_5$ has 8 classes of maximal subgroups:
$X_{5,1} \cong \ZZ_2 \times \Sigma_5$, but $\overline{M}$ of type $X_{5,1}$ implies
$|\overline{G}:\overline{H}| \ge 6 |X_{5}:X_{5,1}|= 6 \cdot 336 = 2016$, a contradiction as in case $X_0$.
$X_{5,2}$ and $X_{5,3}$ are soluble of sizes $2^5 \cdot 3^3$ and $2^8 \cdot 3$,
$X_{5,4} \cong \Alt_6.2^2$ and $X_{5,5} \cong \ZZ_2 \times \PSL_3(2).2$ have bad components,
$X_{5,6} \cong \PSL_3(4).2$ with $2B$-involutions is analyzed below,
$X_{5,7} \cong \PSL_3(4).2$ with $2C$-involutions is analyzed below, but 
$X_{5,8} \cong \PSL_3(4).2$ with $2D$-involutions is out.

$X_{5,6} \cong \PSL_3(4).2$ with $2B$-involutions has 10 classes of maximal subgroups:
$X_{5,6,1} \cong \Sigma_5$, three classes of $\PSL_3(2).2$, three classes of $\Alt_6.2$, soluble groups
of sizes $2^4 \cdot 3^2$ and $2^7 \cdot 3$ and $S \cong \PSL_3(4)$ itself.
Only $X_{5,6,1}$ for $\overline{M}$ remains, but then $|\overline{G}:\overline{H}|\ge 6 \cdot |X_{5,6}:X_{5,6,1}|= 6 \cdot 336 = 2016$
gives a contradiction as before.

$X_{5,7} \cong \PSL_3(4).2$ with $2C$-involutions has 6 classes of maximal subgroups:
a soluble group of size $2^4 \cdot 3^2$, a $\ZZ_2 \times \PSL_3(2)$ , a $\Alt_6.2$ , 
two classes $X_{5,7,4}$ and $X_{5,7,5}$ of shape $2^4:\Sigma_5$ and $S$ itself. 
If $\overline{M}$ is of type $X_{5,7,4}$ or $X_{5,7,5}$, then 
$|\overline{G}:\overline{H}| \ge 6 \cdot |X_{5,7}:X_{5,7,4}| = 6 \cdot 21 = 126$,
but $|\overline{K}| \le 1 +120=121$, as class $2C$ has size 120 in this subgroup.  

$X_7$ with $2C$-involutions has 6 classes of maximal subgroups:
two soluble classes of subgroup sizes $2 \cdot 3^2 \cdot 7$ resp. $2^4 \cdot 3^3$,
two classes $X_{7,3}$ and $X_{7,4}$ of shape $2^4:((3 \times A_5):2)$ and 
two classes containing $\PSL_3(4)$: $\PSL_3(4).3$ and $\PSL_3(4).2 \cong X_{5,7}$. 
Notice, that $\overline{M}$ of type $X_{7,3}$ or $X_{7,4}$ implies 
$3 \in \pi(H)$ by \refOTTprime, but $2C$ involutions invert elements from all classes 
of 3-elements.
\end{bew}

\begin{lemma}
Let $S \cong \PSL_3(q)$, $\PSL_4(q)$, $\PSU_3(q)$, $\PSU_4(q)$, $\Sp_4(q)$, $G_2(q)$, ${}^3D_4(q)$ or ${}^2F_4(q)$
for $q>4$. Then $\overline{G}=\overline{H}$. 
\end{lemma}

\begin{bew}
We use \ref{commuting_in_char_2} as well as $FS_r$-property for $r>2$. Further we use the discussion
of the connected components of the commuting graph $\Gamma_{\cal O}$ together with \ref{simple_group_is_enough}.
Let $x \in H \cap G_0$ be an element of odd prime order $r$, which exists by \refNloopEmbedding. 

In case of $S \cong \PSL_3(q)$, either $r \mid q^2-1$, or $r \mid \frac{q^2+q+1}{(q-1,3)}$. 
If $r \mid q^2-1$, then $\overline{G_0} \subseteq \overline{H}$ by \refsmallCC and
\ref{commuting_in_char_2}. Then $\overline{H}=\overline{G}$ by \ref{simple_group_is_enough}.
If $r \mid \frac{q^2+q+1}{(q-1,3)}$, then $3 \in \pi(H)$ by \ref{H_covers_automizer}, but $3 \mid q^2-1$. 

In case of $S \cong \PSL_4(q)$ the graph $\Gamma_{\cal O}$ is connected by \refsmallCC,
so $\overline{G_0} \subseteq \overline{H}$ by \ref{commuting_in_char_2} and $\overline{H}=\overline{G}$ by \ref{simple_group_is_enough}.

In case of $S \cong \PSU_3(q)$, either $r \mid q^2-1$, or $r \mid \frac{q^2-q+1}{(q+1,3)}$. 
If $r \mid q^2-1$, then $\overline{G_0} \subseteq \overline{H}$ by \refsmallCC and
\ref{commuting_in_char_2}. Then $\overline{H}=\overline{G}$ by \ref{simple_group_is_enough}.
If $r \mid \frac{q^2-q+1}{(q+1,3)}$, then $3 \in \pi(H)$ by \ref{H_covers_automizer}, but $3 \mid q^2-1$. 

In case of $S \cong \PSU_4(q)$ the graph $\Gamma_{\cal O}$ is connected by \refsmallCC,
so $\overline{G_0} \subseteq \overline{H}$ by \ref{commuting_in_char_2} and $\overline{H}=\overline{G}$ by \ref{simple_group_is_enough}.

In case of $S \cong \Sp_4(q)$ either $r \mid q^2-1$ or $r \mid q^2+1$. If $r \mid q^2+1$, then $O_2(\overline{H})=1$: 
No prime divisor of $q^2+1$ divides the order of a parabolic subgroup, so there is no $\{2,r\}$-subgroup of $\overline{G_0}$.
Furthermore the centralizer of an outer involution, which is either $\Sp_4(q^{1/2})$ or ${}^2B_2(q)$, does not contain a torus of size $q^2+1$, which is contained
in $C_{\overline{H}}(\overline{x})$. Therefore we can apply \refHeissPrime. Notice, that neither parabolic subgroups nor subgroups
of type $\Sp_4(q^{1/2})$ or ${}^2B_2(q)$ contain a Sylow-$s$-subgroups for primes $s$ dividing $q+1$. 
Therefore the length of every conjugacy class of involutions is divisible by $s$, so $s \in \pi(H)$, but $s \mid q^2-1$. 
By \refHeissPrime, $s \in \pi(H)$. 
If $r \mid q^2-1$, we have $\overline{G_0} \subseteq \overline{H}$ by \ref{commuting_in_char_2} and \refsmallCC,
so $\overline{H}=\overline{G}$ by \ref{simple_group_is_enough}.

In case of $S \cong G_2(q)$ with $3 \mid q-\varepsilon $ for $\varepsilon \in \{+1,-1\}$, if $r \mid q^2-\varepsilon q+1$
then $3 \in \pi(H)$ by \ref{H_covers_automizer}, and $3 \mid q^2-1$. 
So we  have $r \mid (q^2-1)(q^2+\varepsilon q +1)$ and $\overline{G_0} \subseteq \overline{H}$ by \ref{commuting_in_char_2} and \refsmallCC,
so $\overline{H}=\overline{G}$ by \ref{simple_group_is_enough}.   

In case of $S \cong {}^3D_4(q)$, either $r \mid q^4-q^2+1$ or $r \mid q^6-1$. 
Suppose $r \mid q^4-q^2+1$. Then $\overline{H}$ contains a torus of size $q^4-q^2+1$ from $C_H(x)$, as $C_G(x)$ is soluble.
Therefore $O_2(\overline{H})=1$, so by \refHeissPrime and the list of maximal subgroups of $S$ in \cite{K3D4},
$s \in \pi(H)$ for $s$ some prime divisor of $\frac{q^4+q^2+1}{3}$. 
If $r \mid q^6-1$, then $\overline{G_0} \subseteq \overline{H}$ by \ref{commuting_in_char_2} and \refsmallCC,
so $\overline{H}=\overline{G}$ by \ref{simple_group_is_enough}.

In case of $S \cong {}^2F_4(q)$, either $r \mid q^4-q^2+1$ or $r \mid (q^4-1)(q^3+1)$.
If $r \mid q^4-q^2+1$, then $3 \in \pi(H)$ by \ref{H_covers_automizer}, as $C_H(x)$ contains
the normalizer of a torus either of size $q^2+\sqrt{2q^3}+q+\sqrt{2q}+1$ or $q^2-\sqrt{2q^3}+q-\sqrt{2q}+1$
and $3 \mid q+1$. 
If $r \mid (q^4-1)(q^3+1)$, then $\overline{G_0} \subseteq \overline{H}$ by \ref{commuting_in_char_2} and \refsmallCC,
so $\overline{H}=\overline{G}$ by \ref{simple_group_is_enough}.
\end{bew}

\subsubsection{High rank: the case $q\ge 4$}
In case $q>4$ we use \ref{commuting_in_char_2} as well as $FS_r$-property for $r>2$ by \ref{FS_p_even_char}. 
At this point the discussion of the connected components of the commuting graph $\Gamma_{\cal O}$  becomes essential.
The arguments in case $q=4$ are not that different:
Notice, that we have $FS_r$-property for $r>3$ by \ref{FS_p_even_char}, in particular for $r=5$. 
Recall from \ref{commuting_in_char_2}, that if $FS_3$-property fails, then $5 \in \pi(H)$, so $\overline{H}$
contains a full Sylow-$5$-subgroup of $\overline{G}$ while $q+1=5$ for $q=4$.

\begin{lemma}
Let $S \cong \PSL_n(q)$ or $\PSU_n(q)$ with $n \ge 5$,
    $\Sp_{2n}(q)$ with $n \ge 3$ or
    $\Omega^\pm_{2n}(q)$ for $n \ge 4$ and $q\ge 4$, $q$ even. 
Then $\overline{G}=\overline{H}$. 
\end{lemma}

\begin{bew}
Let $x \in H\cap G_0$ be an element of odd prime order $r$, which exists by \refNloopEmbedding. 

In case of $S \cong \PSL_n(q)$ or $\PSU_n(q)$ we use \refsmallCC.
If $q=4$, then by \refFconnected and \ref{EvenCharTerminals}, we have $\overline{H}=\overline{G}$,
if $5 \in \pi(H)$. 
Then either $\overline{H}=\overline{G}$ by \ref{commuting_in_char_2} and \ref{simple_group_is_enough} or 
$C_{\overline{G}}(\overline{x})$ is a self centralizing torus, on which a prime $p=n$ resp. $p=n-1$ acts. In that case $p$ occurs in $\pi(H)$ by \ref{H_covers_automizer}.
Notice, that in $S$  there is at most one exceptional self centralizing torus, which does not contain
the prime $p$ itself. Therefore elements of order $p$ are in the big connected component and $\overline{H}=\overline{G}$.

In case of $S \cong \Omega^+_{2n}(q)$ or $S \cong \Sp_6(q)$ the graph $\Gamma_{\cal O}$ is connected by \refsmallCC.
If $q=4$, then by \refFconnected and \ref{EvenCharTerminals}, we have $\overline{H}=\overline{G}$,
if $5 \in \pi(H)$. 
So $\overline{H}=\overline{G}$ by \ref{commuting_in_char_2} and \ref{simple_group_is_enough}.

In case of $S \cong \Sp_{2n}(q)$ or $\Omega^-_{2n}(q)$ for $n\ge 4$ we have either
$\overline{H}=\overline{G}$ by \refsmallCC or $n$ is a 2-power and $r \mid q^n+1$. 
If $q=4$, again by \refFconnected and \ref{EvenCharTerminals}, we have $\overline{H}=\overline{G}$,
if $5 \in \pi(H)$. 
We determine the isomorphism type of maximal subgroups $\overline{M}$ of $\overline{G}$, which contain $\overline{H}$. 
Notice, that $\overline{H}_{2'} \ge 4^4+1=257$. 
In the symplectic case we get by \cite{KL} $F^\ast(\overline{M}) \cong \Sp_{n}(q^2)$ or $\Omega^-_{2n}(q)$ or $\overline{M}$ in class ${\cal S}$,
while in the orthogonal case we have $F^\ast(M) \cong \Omega^-_{n}(q^2)$ or in class ${\cal S}$. So in any case $F^\ast(\overline {M})$ 
is a simple group. By \refHsuper, $M$ is a group to a subloop, so we can use \refEnvelopeGroup on $\langle M \cap K \rangle$.
If $F^\ast(\overline{M})$ is a passive group, then $3 \in \pi(H)$ or $5 \in \pi(H)$ with $15 \mid q^4-1$. So $\overline{H}=\overline{G}$
as $\overline{H}$ contains elements of odd order, which are in the big connected component. 
Else $F^\ast(M)$ is a group $\PSL_2(p)$ for a Fermat prime $p\ge 257$ with $p=q^n+1$. The minimal degree of a faithful representation
of $\PSL_2(p)$ is $\frac{p-1}{2}= \frac{q^n}{2}$ by Landazuri-Seitz, but $\overline{G}$ has a faithful module in dimension $2n$. 
As $q\ge 4$, this is absurd.
\end{bew}

\begin{lemma}
Let $S \cong F_4(q)$ $E_6(q)$, ${}^2E_6(q)$, $E_7(q)$ or $E_8(q)$ for $q\ge 4$, $q$ even.
Then $\overline{G}=\overline{H}$. 
\end{lemma}

\begin{bew}
Let $x \in H$ be an element of odd prime order $r$, which exists by \refNloopEmbedding. 

In case of $S \cong F_4(q)$ we use \refsmallCC.
If $q=4$, by \refFconnected and \ref{EvenCharTerminals}, we have $\overline{H}=\overline{G}$,
if $5 \in \pi(H)$.  
Then either $x$ is in the big connected component and $\overline{H}=\overline{G}$ by \ref{commuting_in_char_2} and \ref{simple_group_is_enough}
or $C_{\overline{G}}(\overline{x})$ is a self centralizing torus of size $q^4+1$ or $q^4-q^2+1$. 
The normalizer of the torus of size $q^4-q^2+1$ contains elements of order $3$, as $\overline{G}$ contains 
a subgroup ${}^3D_4(q).3$ with an outer field automorphism of order 3. We can find this field automorphism acting on top of the torus, 
so by \ref{H_covers_automizer}, $3 \in \pi(H)$. As elements of order 3 are in the big connected component, $\overline{H}=\overline{G}$ in this case.
Remains the torus of size $q^4+1$. Let $\overline{M}$ be a maximal subgroup of $\overline{G}$ containing $\overline{H}$. We can use \refEnvelopeGroup on $\langle M \cap K \rangle$ by \refHsuper. As we know, which elements of odd order occure in $\overline{H}$,
there remains only the case, that $\overline{H} \cong \PSL_2(p)$ for $p$ some Fermat prime with $p=q^4+1 \ge 4^4+1=257$.  
By the boundaries of Landazuri-Seitz, a faithful representation of $\PSL_2(p)$ has dimension at least $\frac{p-1}{2} \ge 128$, 
but $S$ has a faithful representation in dimension 26, a contradiction.

In case $S \cong E_6(q)$, ${}^2E_6(q)$, $E_7(q)$ or $E_8(q)$ we use \refsmallCC.
If $q=4$, by \refEconnected and \ref{EvenCharTerminals}, we have $\overline{H}=\overline{G}$,
if $5 \in \pi(H)$. 
Then either $x$ is in the big connected component and $\overline{H}=\overline{G}$ by \ref{commuting_in_char_2} and \ref{simple_group_is_enough}
or $C_{\overline{G}}(\overline{x})$ is a self centralizing torus, on which elements of order $3$ or $5$ act nontrivially.
By \ref{H_covers_automizer} then $3 \in \pi(H)$ or $5 \in \pi(H)$, but elements of order $3$ or $5$ are in the big connected component, so $\overline{H}=\overline{G}$.
\end{bew}

\subsubsection{High rank: the case $q=2$}
The remaining groups are $\PSL_n(2),\PSU_n(2)$ for $n \ge 5$, $\Sp_{2n}(2)$ for $n \ge 3$, $\Omega^\pm_{2n}(2)$ for $n \ge 4$,
 $F_4(2)$, $E_6(2)$, ${}^2E_6(2)$, $E_7(2)$ and $E_8(2)$. 

These cases behaves differently from $q>4$, if $3 \in \pi(H)$ or $5 \in \pi(H)$. 
We use results on centralizer of elements of order 3 and 5, to overcome the
exceptions of \ref{commuting_in_char_2}.

\begin{lemma}
Let $S \cong \PSL_n(2)$ for $n \ge 5$. Then $\overline{H}=\overline{G}$. 
\end{lemma}

\begin{bew}
For $n=5$ we use the list of maximal subgroups and centralizer sizes in \cite{ATLAS}. 
Elements of class $3A$ are terminal by \refFconnected and \ref{EvenCharTerminals}.
Elements of order 7 have soluble centralizer, which contains elements of class 3A.
By the list of maximal subgroup in \cite{ATLAS} this implies $H$ is of type $31:5$,
a contradiction to \refNloopEmbedding, so $n \ge 6$.

Suppose $3 \in \pi(H)$. Let $V$ be the natural $n$-dimensional $\GF(2)$-module for $S$. 
By \refFconnected and \ref{EvenCharTerminals} there exists a terminal element $\overline{t}$ of order 3,
which has $\dim [V,\overline{t}]=2$. Therefore $\overline{H}$ does not contain a Sylow-3-subgroup of $G$, so $FS_3$-property fails.
How can $FS_3$-property fail? From the structure of centralizers of semisimple elements and the structure of nonsoluble subloops
we conclude as in \ref{even_char_centralizers}, that some $\overline{y} \in \overline{H}, o(\overline{y})=3$ exists,
such that $C_{\overline{G}}(\overline{y})$ contains a subnormal subgroup isomorphic to $\Sigma_5$. 
By \ref{GL_n_2-35Centralizer}, $C_{\overline{G}}(\overline{y}) \cong \GL_{m/2}(4) \times \SL_{n-m}(2)$ with $m:=\dim[V,y]$ 
for $V$ the natural $n$-dimensional $\GF(2)$-module of $G$. In particular components of type $\Alt_5 \cong \PSL_2(4)$ occure only for $m=4$. 
If $m \ge 6$, $\overline{H}$ covers the $\SL_{n-m}(2)$ acting on $C_V(\overline{y})$ by \refpassivecentralizingcomponents. We conclude, that then $\overline{H}$ contains
elements, which are conjugate to $\overline{t}$, so $\overline{H}=\overline{G}$ for $n \ge 6$ and $3 \in \pi(H)$. 

So let $5 \in \pi(H)$, $\overline{x} \in \overline{H}$ of order 5 and consider
the action of $C_{\overline{G}}(\overline{x})$. Let $m:=\dim [V,\overline{x}]$. By \ref{GL_n_2-35Centralizer}, 
$C_{\overline{G}}(\overline{x}) \cong \GL_{m/4}(16) \times \SL_{n-m}(2)$. 
We conclude, that $O_3(C_{\overline{G}}(\overline{x})) \ne 1$, so $3 \in \pi(H)$ by \refOTTprime. 

Finally let $3 \not \in \pi(H) \not\ni 5$. From \refsmallCC we conclude, that either $n$ or $n-1$ is a prime
and $\overline{H}$ contains a torus of size $2^n-1$ resp. $2^{n-1}-1$. Then by \ref{H_covers_automizer} either $n \in \pi(H)$
or $n-1 \in \pi(H)$, and $n$ resp. $n-1$ are in the connected component containing all elements of order $3$ and $5$.
We then get a contradiction, as this implies $3 \in \pi(H)$ or $5 \in \pi(H)$. 
\end{bew}

\begin{lemma}
Let $S \cong \Sp_n(2)$ for $n \ge 6$. Then $\overline{H}=\overline{G}$.
\end{lemma}

\begin{bew}
Let $V$ be the natural $n$-dimensional module of $\overline{G}$. Recall $\Out(G)=1$. 

Suppose $3 \in \pi(H)$. By \refFconnected there exists a  terminal element $\overline{t}$ of order 3,
with $\dim [V,\overline{t}]=2$. Therefore $\overline{H}$ does not contain a Sylow-3-subgroup of $G$, so $FS_3$-property fails.
We conclude as in \ref{even_char_centralizers}, that some $\overline{y} \in \overline{H}, o(\overline{y})=3$ exists,
such that $C_{\overline{G}}(\overline{y})$ contains a subnormal subgroup isomorphic to $\Sigma_5$. 
Let $\overline{y} \in \overline{H}$ with $o(\overline{y})=3$. 
By \ref{Sp_n_2-35Centralizer}, $C_{\overline{G}}(\overline{z}) \cong \GU_{m/2}(2) \times \Sp_{n-m}(2)$,
so no such subnormal subgroup occurs and $FS_3$-property holds, so $3 \not\in \pi(H)$.

Suppose $5 \in \pi(H)$. Let $\overline{x} \in \overline{H}$ be of order 5 and $m:=\dim [V,x]$. 
By \ref{Sp_n_2-35Centralizer}, $C_{\overline{G}}(\overline{y}) \cong \GU_{m/4}(4) \times \Sp_{n-m}(2)$.
If $\dim C_V(x) >0$, $\overline{H}$ covers the $\Sp_{n-m}(2)$-factor of $C_{\overline{G}}(\overline{y})$ by \refpassivecentralizingcomponents,
so $3 \in \pi(H)$. If $m \ge 12$, then $\overline{H}$ covers the $\GU_{m/4}(4)$-factor too and $3 \in \pi(H)$.
So $n=m=8$, but then the centralizer has structure $\ZZ_5 \times \Alt_5$ and no subgroup $\Sigma_5$ occurs in this centralizer,
so $3 \in \pi(H)$ in this case too.   

If now $3 \not\in \pi(H) \not\ni 5$, we use \refsmallCC for the connected
components of $\Gamma_{\cal O}$, which do not contain elements of order 3 or 5.
Let $\overline{x} \in \overline{H}$ be an element of odd order. We conclude, that either
$n$ is 2-power and $o(x)$ divides $q^{n/2}+1$ or $n=2p$ for a prime $p$ and 
$n$ divides $2^p-1$. In this last case $p \in \pi(H)$ by \ref{H_covers_automizer} and $p$ is in the big connected
component, so $\overline{H}=\overline{G}$.
Let $\overline{M}$ be a maximal subgroup containing $\overline{H}$, so $M$ is a group to a subloop by \refHsuper.
We get by \cite{KL}, that $F^\ast(\overline{M}) \cong \Sp_{n/2}(q^2)$ or $\Omega^-_{n}(q)$ or $\overline{M}$ in class ${\cal S}$.
In any case $F^\ast(\overline{M})$ is a simple group.  If this group is passive, either $3 \in \pi(H)$ or $5 \in \pi(H)$ and $\overline{H}=\overline{G}$.
Else $\overline{M} \cong \PGL_2(p)$ for a Fermat prime $p \ge 17$. We have $p \mid 2^{n/2}+1$, $G$ has an $n$-dimensional $\GF(2)$-module,
but the minimal representation degree of $\PGL_2(p)$ is $p-1$. Now $2^{n/2} \le n$, so $n \le 4$, a contradiction. 
Notice, that the group $\Sp_8(2)$ actually has a maximal subgroup isomorphic to $\PSL_2(17)$, but no $\PGL_2(17)$. 
\end{bew}

\begin{lemma}
Let $S \cong \Omega^+_n(2)$ for $n \ge 8$. Then $\overline{H}=\overline{G}$.
\end{lemma}

\begin{bew}
Let $V$ be the natural $n$-dimensional module of $\overline{G}$. Recall, that $\Out(G)=\ZZ_2$ for $n \ge 10$
and $\Sigma_3$ for $n=8$. 

Suppose $3 \in \pi(H)$. By \refFconnected and \ref{EvenCharTerminals} there exists a terminal element $\overline{t}$ of order 3,
which have $\dim [V,\overline{t}]=2$. Therefore $\overline{H}$ does not contain a Sylow-3-subgroup of $G$, so $FS_3$-property fails.
We conclude as in \ref{even_char_centralizers}, that some $\overline{y} \in \overline{H}, o(\overline{y})=3$ exists,
such that $C_{\overline{G}}(\overline{y})$ contains a subnormal subgroup isomorphic to $\Sigma_5$. 
Let $\overline{y} \in \overline{H}$ with $o(\overline{y})=3$. 
By \ref{OP_n_2-35Centralizer}, $O^2(C_{\overline{G}}(\overline{y})) \cong (\GU_{m/2}(2))' \times \Omega^{\varepsilon_1}_{n-m}(2)$ for $\varepsilon_1=(-1)^{m/2}$,
if $\overline{y}$ is an element of $\Omega^+_n(2)$. If $\overline{y}$ is outside of $\Omega^+_8(2)$, we
use \cite{ATLAS} for the structure of $O^2(C_{\overline{G}}(\overline{y}))$. 
In any case a subnormal $\Alt_5$ exists only for $n-m=4$, with $\varepsilon_1=-1$. 
In that case $\overline{H}$ covers the subgroup $\GU_{m/2}(2)$ by \refpassivecentralizingcomponents, 
so $\overline{H}$ contains an element, which is conjugate to $\overline{x}$ and $\overline{H}=\overline{G}$. 
Therefore either $\overline{H}=\overline{G}$ or $FS_3$-property holds, so $3 \not\in \pi(H)$.

Suppose $5 \in \pi(H)$. Let $\overline{x} \in \overline{H}$ be of order 5 and $m:=\dim [V,\overline{x}]$. 
By \ref{OP_n_2-35Centralizer}, $C_{\overline{G}}(\overline{y}) \cong \GU_{m/4}(4) \times \Omega^{\varepsilon_2}_{n-m}(2)$
for $\varepsilon_2=(-1)^{m/4}$.
Suppose $n-m \ge 6$. Then $\Omega^{\varepsilon_2}_{n-m}(2)$ is passive, so by \refpassivecentralizingcomponents,
$3 \in \pi(H)$. 
If $(n-m,\varepsilon_2)= (4,+1) $ or $(2,-1)$, for the same reason $3 \in \pi(H)$. 
If $m \ge 12$, then $\overline{H}$ covers the $\GU_{m/4}(4)$-factor too  by \refpassivecentralizingcomponents and $3 \in \pi(H)$.
So $m \le 8$ and $(n-m,\varepsilon_2) \in \{ (0,+1)=(0,-1), (2,+1), (4,-1) \}$. This gives the groups
$O_8^+(2)$ and $O^+_{10}(2)$.  
In both cases $\overline{H}$ contains the normalizer of a Sylow-5-subgroup. 
In case of $S \cong \Omega^+_{10}(2)$ we get $3 \in \pi(H)$: We check the list of maximal subgroup in \cite{ATLAS}
and use \refHsuper with \refEnvelopeGroup, to get $3 \in \pi(H)$. 
If $S \cong \Omega^+_8(2)$,  $3 \in \pi(H)$ by \ref{H_covers_automizer},
if $\overline{G} \cong \Omega^+_8(2).\Sigma_3$, so $|\overline{G}:\overline{G_0}| \le 2$.
Calculation of structure constants within $\Omega^+_8(2).2$ reveals, that nontrivial elements of $\overline{K}$ are
either in class $2A$ or $2F$ (Notation as in \cite{ATLAS}), as the other classes of involutions invert elements of order 5. 
Therefore $|\overline{K}| \le 1 + 1575 + 120 = 1796$. On the other hand $|\overline{G}:\overline{H}| \ge 2 \cdot 3^5 \cdot 7= 3402$, a contradiction. 
 
If now $3 \not\in \pi(H) \not\ni 5$, we use \refsmallCC for the connected
components of $\Gamma_{\cal O}$, which do not contain elements of order 3 or 5.
In particular there exists a prime $p$ with $n=2p$ or $n=2p+2$ and the connected component contains
elements of prime order $r$ for all prime divisors $r$ of $2^p-1$. By \ref{H_covers_automizer}, $r \in \pi(H)$
implies $p \in \pi(H)$, while $p$ itself is in the big connected component. Therefore $\overline{H}=\overline{G}$. 
\end{bew}

\begin{lemma}
Let $S \cong \Omega^-_n(2)$ for $n \ge 8$. Then $\overline{H}=\overline{G}$.
\end{lemma}

\begin{bew}
Let $V$ be the natural $n$-dimensional module of $\overline{G}$. Recall, that $\Out(G)=\ZZ_2$.

Suppose $3 \in \pi(H)$. By \refFconnected and \ref{EvenCharTerminals} there exists terminal elements $\overline{t}$ of order 3,
with $\dim [V,\overline{t}]=2$. Therefore $\overline{H}$ does not contain a Sylow-3-subgroup of $G$, so $FS_3$-property fails.
We conclude as in \ref{even_char_centralizers}, that some $\overline{y} \in \overline{H}, o(\overline{y})=3$ exists,
such that $C_{\overline{G}}(\overline{y})$ contains a subnormal subgroup isomorphic to $\Sigma_5$. 
Let $\overline{y} \in \overline{H}$ with $o(\overline{z})=3$. 
By \ref{OM_n_2-35Centralizer}, $O^2(C_{\overline{G}}(\overline{y})) \cong (\GU_{m/2}(2))' \times \Omega^{\varepsilon_1}_{n-m}(2)$ for $\varepsilon_1=(-1)^{1+m/2}$.
Therefore a subnormal $\Alt_5$ exists only for $n-m=4$, with $\varepsilon_1=-1$. 
In that case $\overline{H}$ covers the subgroup $\GU_{m/2}(2)$ by \refpassivecentralizingcomponents, 
so $\overline{H}$ contains an element, which is conjugate to $\overline{x}$ and $\overline{H}=\overline{G}$. 
Therefore either $\overline{H}=\overline{G}$ or $FS_3$-property holds, so $3 \not\in \pi(H)$.

Suppose $5 \in \pi(H)$. Let $\overline{x} \in \overline{H}$ be of order 5 and $m:=\dim [V,x]$. 
By \ref{OM_n_2-35Centralizer}, $C_{\overline{G}}(\overline{y}) \cong \GU_{m/4}(4) \times \Omega^{\varepsilon_2}_{n-m}(2)$
for $\varepsilon_2=(-1)^{1+m/4}$.
Suppose $n-m \ge 6$. Then $\Omega^{\varepsilon_2}_{n-m}(2)$ is passive, so by \refpassivecentralizingcomponents,
$3 \in \pi(H)$. 
If $(n-m,\varepsilon_2)= (4,+1) $ or $(2,-1)$, for the same reason $3 \in \pi(H)$. 
If $m \ge 12$, then $\overline{H}$ covers the $\GU_{m/4}(4)$-factor too  by \refpassivecentralizingcomponents and $3 \in \pi(H)$.
The case $m=4$ gives a contradiction as then $C_V(\overline{x})$ is a $O^+_4(2)$-space.
In case $m=8$, $[V,\overline{x}]$ is a $O^+_8(2)$-space. Therefore $C_V(\overline{x})$ has to be an
$O_4^-(2)$-space and $S \cong \Omega^-_{12}(2)$. But from the centralizer structure we conclude, that $FS_5$-property holds.
So $\overline{H}$ contains a Sylow-5-subgroup of $\overline{G}$ and there are other elements of order 5 in $\overline{H}$
which imply $3 \in \pi(H)$. 

If now $3 \not\in \pi(H) \not\ni 5$, we use \refsmallCC for the connected components of $\Gamma_{\cal O}$, which do not contain elements of order 3 or 5.
Then either $n$ or $n-2$ is a 2-power and the connected component contains elements of order $r$
for primes $r$ dividing $2^{n/2}+1$ resp. $2^{n/2-1}+1$. 
Let $\overline{M}$ be a maximal subgroup of $\overline{G}$ containing $\overline{H}$.
By \cite{KL} we get $\overline{M}$ of type $O^-_{n/2}(4)$ or in class ${\cal S}$, if $n$ is a 2-power.
If $n-2$ is a 2-power, then $\overline{M}$ is of type $\Sp_{2n-2}(2)$, a parabolic of type $2^{n-2}:O^-_{n-2}(2)$ 
or $\overline{M}$ in class ${\cal S}$.  
Using \refHsuper and \refEnvelopeGroup we get $3 \in \pi(H)$, $5 \in \pi(H)$ or $\overline{M} \cong \PGL_2(p)$
for a Fermat prime $p$. Then $p=2^{n/2}+1$ or $p=2^{n/2-1}+1$ and $\PGL_2(p)$ has a faithful representation in
degree at least $p-1$, but $\overline{G}$ has an $n$-dimensional module, so $2^{n/2-1} \le n$, which gives
contradictions: either $p=9$ or $n \ge 10$. 
\end{bew}

\begin{lemma}
Let $S \cong \PSU_n(2)$ for $n \ge 5$. Then $\overline{H}=\overline{G}$.
\end{lemma}

\begin{bew}
Recall, that $\Out(S) \cong \ZZ_2$ or $\Sigma_3$ depending on whether $n$ is divisible by 3.

Suppose $3 \in \pi(H)$. By \refFconnected and \ref{EvenCharTerminals} there are terminal elements of order $3$, so 
$\overline{H}$ does not contain a Sylow-3-subgroup of $\overline{G}$, so $FS_3$-property fails.
So there exists some $x \in \overline{H}$, $o(\overline{x})=3$, such that $C_{\overline{G}}(\overline{x})$ contains
a subnormal $\Sigma_5\cong \PGaL_2(4)$. 
We use \ref{unitary_3_centralizers} for the description of centralizers of elements of order 3. In particular
we see, that no subnormal $\Sigma_5$ exists, so $FS_3$-property holds and $\overline{H}=\overline{G}$. 

Suppose $5 \in \pi(H)$. By \ref{unitary_5_centralizers} and \ref{FS_p_criterion}, property $FS_5$ holds,
so $\overline{H}$ contains a Sylow-5-subgroup of $\overline{G}$. 
In particular $\overline{H}$ contains an element $\overline{x}$ of order $5$, such that
$\dim [V,\tilde{x}]=4$, for $\tilde{x}$ some preimage of $\overline{x}$ in $\GU_n(2)$ and $V$ the natural $\GF(4) \GU_n(2)$-module.
For this element, $O_3(C_{\overline{G}}(\overline{x})) \ne 1$, so $3 \in \pi(H)$ by \refOTTprime. 

So suppose $3 \not\in \pi(H) \not\ni 5$. We use \refsmallCC for the connected components of $\Gamma_O$, which do not contain
elements of order 3 or 5.
Let $\overline{x} \in \overline{H}$, $o(\overline{x})$ some odd prime $r$. 
By \ref{commuting_in_char_2} we conclude, that a prime $p$ exists with either $n=p$ or $n-1=p$ and
$r \mid \frac{2^p+1}{3}$. Then $\overline{H}$ contains a torus of size $\frac{2^p+1}{3}$, on which a subgroup of size $p$ acts.
By \ref{H_covers_automizer} then $p \in \pi(H)$, but $p$ is in the connected component containing elements of order 3 and 5,
so $3 \in \pi(H)$ or $5 \in \pi(H)$. 
\end{bew}

\begin{lemma}
Let $S \cong F_4(2)$, $E_6(2)$ or ${}^2E_6(2)$. Then $\overline{H}=\overline{G}$. 
\end{lemma}

\begin{bew}
Case $S\cong F_4(2)$: We use \cite{ATLAS}. By inspection of centralizers of 3-elements we get $FS_3$-property,
so together with \refFconnected and \ref{EvenCharTerminals} we get $\overline{H}=\overline{G}$, if $3 \in \pi(H)$.
If $5 \in \pi(H)$ we get $3 \in \pi(H)$: The centralizer of a $5A$-element has structure $\ZZ_5 \times \Sp_4(2)$, so we get elements of order 3 into $\overline{H}$.
By \refTNloopEnvelope there remains only the case $17 \in \pi(H)$. Subgroups of order 17 are self centralizing.
Since we cannot exclude the existence of a $\PGL_2(17)$, we count involutions. We get $|\overline{K}| \le 96648112$,
as elements of class $2C$ invert elements of order 17. On the other hand $\overline{H} \cong \ZZ_{17}:\ZZ_{16}$,
so $|\overline{G}:\overline{H}| \ge 2^{19}\cdot 3^6 \cdot 5^2 \cdot 7^2 \cdot 13 = 6086629785600$.

Case $S \cong {}^2E_6(2)$. We use the character table of ${}^2E_6(2).3$ and ${}^2E_6(2).2$ provided by GAP. 
We first establish the $FS_3$-property:  
Let $y \in X={}^2E_6(2).3$ be an element of order 3. We claim, that $C_X(y)$ does not contain a subnormal 
$\Alt_5 \cong \PSL_2(4)$. We show this, using the list of conjugacy classes and centralizer sizes of $X$, which we get from the character table: 
If $C_G(y)$ contains no elements of order 5 or elements of order 11 or 19, the statement is obvious:
Elements of order 11 or 19 do not commute with elements of order 5, but cannot permute components of type $\Alt_5$ nontrivially.
Remains only one class of elements (class 3B), with corresponding centralizer containing a subgroup of order $3^8$.  As elements of order $5$
commute in $G$ with 3-groups of size at most $3^3$, we have $FS_3$-property. 
By \refEconnected and  \ref{EvenCharTerminals}, we have $\overline{H}=\overline{G}$, if $3 \in \pi(H)$. 
If $5 \in \pi(H)$, the centralizer of a $5A$-element is $\ZZ_5 \times \Alt_8$, which is contained in a maximal subgroup $\Omega^-_{10}(2)$. 
Therefore $5 \in \pi(H)$ implies $3 \in \pi(H)$. 
By \refNloopEmbedding there remains only the case $17 \in \pi(H)$. 
Subgroups of order 17 are self centralizing in ${}^2E_6(2)$, but not in ${}^2E_6(2).3$,
so $|\overline{G}:\overline{G_0}| \le 2$. 
Let $\overline{M}$ be a maximal subgroup of $\overline{G}$ containing $\overline{H}$. If $O_2(\overline{M}) \ne 1$,
we have either  $\overline{M}$ parabolic or the centralizer of an outer involution.
If $\overline{M}$ is a  parabolic subgroup, then $3 \in \pi(H)$ by \refHsuper and \refEnvelopeGroup.
Of the two outer classes of involutions only the class with centralizer $F_4(2)$ does not invert elements of order 17.
So $\overline{M} \cong F_4(2)$ and $3 \in \pi(H)$ by \refHsuper and \refEnvelopeGroup.
So $O_2(\overline{M})=1$ and we can use \refHeissPrime.
This shows $3 \in \pi(H)$ and $19 \in \pi(H)$, so $\overline{H}=\overline{G}$. 

Case $S \cong E_6(2)$. Recall $\Out(S) \cong \ZZ_2$. 
We show $FS_3$-property, using the character table of $E_6(2)$ as provided in GAP.
Elements of class $3A$ and $3B$ commute with elements of order 31 resp. 17, while elements of of order 31 and 17
do not commute with elements of order 5, so the centralizers of $3A$ and $3B$-elements do not contain subnormal
$\Alt_5=\PSL_2(4)$-subgroups.
Remains the centralizer of a $3C$-Element, which has size $2^9 \cdot 3^6 \cdot 5 \cdot 7$. As the centralizer of
a $5A$-element has size $2^6 \cdot 3^4 \cdot 5^2$, we get a contradiction: A component of type $\PSL_2(5)$ or $\SL_2(5)$
would be normal in $C_S(3C)$, so a Sylow-3-subgroup of size $3^6$ acts on it. 
This gives a contradiction, as the kernel of this action has size at most $3^4$.
By \refEconnected and  \ref{EvenCharTerminals}, we have $\overline{H}=\overline{G}$, if $3 \in \pi(H)$.  
If $5 \in \pi(H)$, then $3 \in \pi(H)$. From the existence of a Levi complement $\Omega^+_{10}(2)$, 
\ref{OP_n_2-35Centralizer} and the centralizer size we conclude, that $C_S(5A) \cong \ZZ_5 \times U_4(2)$, so $5 \in \pi(H)$ implies $3 \in \pi(H)$.
By \refNloopEmbedding there remains only the case $17 \in \pi(H)$, but this time $17 \in \pi(H)$ implies $3 \in \pi(H)$ by the centralizer size of $3 \cdot 17$.
\end{bew}

\begin{lemma}
Let $S \cong E_7(2)$ or $E_8(2)$. Then $\overline{H}=\overline{G}$. 
\end{lemma}

\begin{bew}
Recall $\Out(S) =1$ and suppose $O_2(\overline{H}) \ne 1$. 
Then $\overline{H}$ is contained in a maximal parabolic, so by \refHsuper, some maximal parabolic $\overline{P}$
is a group to a subloop. We use \refEnvelopeGroup on it, to get this subloop soluble. In particular
$\overline{H}$ contains a Sylow-3-subgroup of $\overline{P}$. 
By \refEconnected we get a connected conjugacy class of elements of order 3,
which has a centralizer of type $\ZZ_3 \times \Omega^+_{12}(2)$ resp. $\ZZ_3 \times E_7(2)$,
so this conjugacy class is terminal by \ref{EvenCharTerminals}. But any maximal parabolic subgroup contains
such elements, as these elements come from a $\PSL_2(2)$, which is generated by root subgroups $X_\alpha,X_{-\alpha}$.
Therefore $\overline{H}=\overline{G}$, if $O_2(\overline{H}) \ne 1$. 

So $O_2(\overline{H}) =1$ and we would like to use \refHeissPrime.
As centralizers of involutions are contained in maximal parabolics, we see from \cite{ATLAS}, 
p. 219 and p.235, that $3$ is a prime, for which we can apply \refHeissPrime. So $\overline{H}$
contains a Sylow-3-subgroup of $\overline{G}$. 
Since we showed above, that $\overline{G}$ contains terminal elements of order 3, $\overline{H}=\overline{G}$.
\end{bew}

\end{document}